\documentclass[twocolumn]{autart}  
\usepackage{amsmath}
\usepackage{amssymb}
\usepackage{picins}
\usepackage{xcolor}

\usepackage{graphics} 
\usepackage{epsfig} 
\usepackage{epstopdf}
\usepackage{mathptmx} 
\usepackage{times} 
\usepackage{amsmath} 
\usepackage{amssymb}  
\usepackage{breqn}
\usepackage{bbm}
\usepackage{dsfont}

\newtheorem{ex}{Example}

\newcommand{\mcl}[1]{\mathcal{ #1}}
\newcommand{\mbb}[1]{\mathbb{ #1}}
\newcommand{\mbf}[1]{\mathbf{ #1}}

\newcommand\Item[1][]{%
	\ifx\relax#1\relax  \item \else \item[#1] \fi
	\abovedisplayskip=0pt\abovedisplayshortskip=0pt~\vspace*{-\baselineskip}}

\newcommand{\R}{\mathbb{R}}

\newcommand{\N}{\mathbb{N}}

\newcommand{\eps}{\varepsilon}

\begin{document}
		
		\begin{frontmatter}
	\title{A Generalization of Bellman's Equation with Application to Path Planning, Obstacle Avoidance and Invariant Set Estimation}

	\author[Rome]{Morgan Jones}\ead{morgan.c.jones@asu.edu},    
	\author[Rome]{Matthew M.Peet}\ead{mpeet@asu.edu},               
	
	\address[Rome]{\textsc{SEMTE at Arizona State University, Tempe, AZ, 85287-6106 USA }}             


		\begin{abstract}                          
The standard Dynamic Programming (DP) formulation can be used to solve Multi-Stage Optimization Problems (MSOP's) with additively separable objective functions. In this paper we consider a larger class of MSOP's with monotonically backward separable objective functions; additively separable functions being a special case of monotonically backward separable functions. We propose a necessary and sufficient condition, utilizing a generalization of Bellman's equation, for a solution of a MSOP, with a monotonically backward separable cost function, to be optimal. Moreover, we show that this proposed condition can be used to efficiently compute optimal solutions for two important MSOP's; the optimal path for Dubin's car with obstacle avoidance, and the maximal invariant set for discrete time systems.
		\end{abstract}
		
	\end{frontmatter}

\section{Introduction}

Throughout Engineering, Economics, and Mathematics many problems can be formulated as Multi-Stage Optimization Problems (MSOP's):
	 \begin{align*}
 \min &\bigg\{  J( u(0),...,u(T-1) , x(0),...,x(T)) \bigg\} \\ \nonumber
& x(0)=x_0, \text{ } x(t+1)=f(x(t),u(t),t) \text{ for  } t={0},..,T-1 \\ \nonumber
&   x(t) \in X_t \subset \mathbb{R}^n, \text{ } u(t) \in U \subset \mathbb{R}^m \text{ for } t={0},..,T. 
\end{align*}
 Such problems consist of 1) a cost function $J: \R^{m \times T} \times \R^{n \times (T+1)} \to \R$, 2) an underlying discrete-time dynamical system governed by the plant equation $f: \R^n \times \R^m \times \N \to \R^n$, 3) a state space $X_t \subset \R^n$, 4) an admissible input space $U \subset \R^m$, and 5) a terminal time $T>0$. Examples of such optimization problems include: optimal battery scheduling to minimize consumer electricity bills \cite{jones2017solving}; energy-optimal speed planning for road vehicles \cite{zeng2018globally}; optimal maintenance of manufacturing systems \cite{liu2019manufacturing}; etc.

MSOP's are members of the class of constrained nonlinear optimization problems. Such optimization problems can be solved using nonlinear solvers such as SNOPT \cite{gill2005snopt} over small time horizons. However, the most commonly used class of methods for solving MSOP's is Dynamic Programming (DP) \cite{bertsekas1995dynamic}. DP methods exploit the structure of MSOP's to decompose the optimization problem into lower dimensional sub-problems that can be solved recursively to give the solution to the original higher dimensional MSOP. Typically, DP is used to solve problems with cost functions of the form $J(\mbf u, \mbf x)= \sum_{t=0}^{T-1} c_t(x(t),u(t)) + c_T(x(T))$. These functions (Defn.~\ref{defn: additively seperbale function}) are called additively separable functions, as they can be additively separated into sub-functions, each of which only depend on a single time-stage, $t \in \{0,...,T\}$. In the additively separable case it was shown in \cite{bellman1966dynamic} that if we can find a function $F$ that satisfies Bellman's Equation,
\begin{align*}
F(x,T) & = c_T(x) \quad \text{for all } x \in X_T \\ \nonumber
F(x,t) & = \inf_{u \in \Gamma_{x,t}} \bigg\{ c_t(x,u)+F(f(x,u,t),t+1) \bigg\} \\ \nonumber
& \qquad \qquad \qquad \text{for all } x \in X_t, t \in \{0,..,T-1\},
\end{align*}
where $\Gamma_{x,t}:= \{u\in U: f(x,u,t) \in X_t\}$, then a necessary and sufficient condition for a feasible input and state sequence, $\mbf u=(u(0),...,u(T-1))$ and $\mbf x =(x(0),...,x(T))$, to be optimal is 
\begin{align*}
 u(t) \in \arg \inf_{ u \in \Gamma_{x(t),t} } \bigg\{ c_t(x(t),u) & + F(f(x(t),u,t),t+1)  \bigg\}\\
& \text{for all }  t \in \{0,..,T-1\}.
\end{align*} We consider MSOP's with cost functions of the more general form $J(\mbf u, \mbf x)= \phi_0(x(0),u(0), \phi_1(x(1),u(1), \dots \phi_T(x(T)) \dots ))$, where maps $\phi_t:X_t \times U \times \R \to \R $ are monotonic in their third argument for $t=0,\cdots T-1$. Such functions are called monotonically backward separable, defined in Definition \ref{defn: montonic backward seperbale function}, and shown to contain the class of additively separable functions in Lemma \ref{lem: add sep is back sep}. For MSOP's with monotonically backward separable cost functions we show in Theorem \ref{thm: GBE} that if we can find a function $V$ that satisfies
\begin{align}\label{intro:GBE}
V(x,T) & = \phi_T(x) \quad \text{for all } x \in X_T \\ \nonumber
V(x,t) & = \inf_{u \in \Gamma_{x,t}} \bigg\{ \phi_t(x,u,V(f(x,u,t),t+1)) \bigg\} \\ \nonumber
& \qquad \qquad \qquad \qquad \text{for all } x \in X_t, t \in \{0,..,T-1\},
\end{align}
where $\Gamma_{x,t}:= \{u\in U: f(x,u,t) \in X_t\}$, then a necessary and sufficient for a feasible input and state sequence, $\mbf u=(u(0),...,u(T-1))$ and $\mbf x =(x(0),...,x(T))$, to be optimal is
\begin{align*}
& u(t) \in \arg \inf_{ u \in \Gamma_{x(t),t} } \bigg\{ \phi_t\left(x(t),  u,V(f(x(t),u,t),t+1)\right)  \bigg\}\\
& \hspace{4.5cm} \text{for all }  t \in \{0,..,T-1\}.
\end{align*}
Equation~\eqref{intro:GBE} can be thought of as a generalization of Bellman's Equation; as it is shown in Corollary \ref{cor: GBE implies BE} that in the special case when the cost function is additively separable Equation~\eqref{intro:GBE} reduces to Bellman's Equation. We therefore refer to Equation~\eqref{intro:GBE} as the Generalized Bellman's Equation (GBE). Through several examples we show a solution, $V$, to the GBE can be obtained numerically by recursively solving the GBE backwards in time for each element of $X_t$, the same way Bellman's Equation is solved, thereby extending traditional DP methods to solve a larger class of MSOP's with non-additively separable cost functions. Moreover, in Section \ref{sec: compare} it is shown how Approximate Dynamic Programming (ADP) methods can be modified to solve the GBE.

By recursively solving the GBE~\eqref{intro:GBE} it is possible to synthesize optimal input sequences for many important practical problems. In this paper we consider two such problems; path planning with obstacle avoidance and maximal invariant sets. First, we define the path planning problem as the search for a sequence of inputs that drives a dynamical system to a target set in minimum time while avoiding obstacles defined by subsets of the state-space. In Section~\ref{sec: path planning} we show that such problems can be formulated as an MSOP with monotonically backward separable objective, of form $J(\mbf u, \mbf x)= \min \left\{ \inf \left\{t \in [0,T] : x(t) \in S \right\}, T \right\}$, implying that the solution to the path planning problem can be found using the solution to the GBE. Similarly, in Section~\ref{sec: invariant sets} we show that computation of maximal invariant sets can be formulated as an MSOP with monotonically backward separable objective of form $J(\mbf u, \mbf x)=\max\{\max_{0 \le k \le T-1}\{c_k(u(k),x(k))\}, c_T(x(T))\}$.

Path planning with obstacle avoidance has been extensively studied (see surveys \cite{dreyfus1969appraisal} \cite{gallo1988shortest}) and has many applications; including UAV surveillance \cite{xie2019optimal}. In \cite{rippel2005fast} the path planning problem is separated into two separate problems: the ``geometric problem", in which the shortest curve, $\tilde{x}(t)$, between the initial set and target set is calculated, and the ``tracking problem", in which a controller, $u(t)$, is synthesized so that $\sum_{t=0}^T ||x(t) - \tilde{x}(t)||_2^2$ is minimized, where $x(t+1)=f(x(t),u(t),t)$ and $|| \cdot||_2$ is the Euclidean norm. Separating the path planning problem allows for the use of efficient algorithms such as $A^*$-search or tangent graphs \cite{liu1992path} to solve the ``geometric problem" and LQR control to solve the ``tracking problem", however, there is no guaranteed that this method will produce the true solution to the original path planning problem. The same approach is used in \cite{cowlagi2011hierarchical}, where it is shown through numerical examples that a controller closer to optimality can be derived when the state space is augmented with historic trajectory information. Our approach of using the GBE to solve the path planning does not separate the problem into the ``geometric or ``tracking" problem and thus does not require any state augmentation. For systems described in continuous time (rather than the discrete systems considered in this paper) with obstacles that satisfy certain boundary curvature assumptions, assumptions not made in this paper, it has been shown in \cite{savkin2013reactive} that a path planning sliding mode controller can be efficiently computed. Furthermore, this sliding mode controller can be used for effective path planning in unknown environments, a case not considered in this paper.

The GBE can also be used in the application of computing the Finite Time Horizon Maximal Invariant Set (FTHMIS), defined as the largest set of initial conditions for a discrete time process such that there exists a feasible input sequence for which the state of the system never violates a time-varying constraint. Knowledge of this set can be used to design controllers that ensure the system never violates given safety constraints.  We show that FTHMIS's are equivalent to the sublevel set of solutions to the GBE. To the best of the authors knowledge the problem of computing FTHMIS's has not previously been addressed in the literature. However, a proposed methodology for computing maximal invariant sets over infinite time horizons can be found in \cite{xue2018robust,esterhuizen2019maximal,wang2019computation}. Similar continuous-time formulations of this problem can be found in \cite{jones2019using,jones2019relaxing}. 

Substantial work on generalizations of Bellman's Equation for both infinite and finite time MSOP's can be found in~\cite{bertsekas2018abstract}. Our work differs from \cite{bertsekas2018abstract} as rather than attempting to generalize the "Bellman's operator``, as \cite{bertsekas2018abstract} does, we consider a wider class of cost functions associated with MSOP's, introducing monotonically backward separable cost functions, leading to a derivation of the GBE~\eqref{intro:GBE}. Unlike in~\cite{bertsekas2018abstract}, we formalize the link between the cost function of an MSOP and the GBE~\eqref{intro:GBE}. Other examples in the literature of MSOP's with non-additively separable cost functions can be found in the pioneering work of Li~\cite{li1991extension,li1990new,li1990multilevel,li1990multiple}. Li considered MSOP's with $k$-separable cost functions; functions of the form $J(\mbf u, \mbf x)=H(J_1(\mbf u, \mbf x),...,J_k(\mbf u, \mbf x))$, where $H:\R^k \to \R$ is strictly increasing and differentiable, and each of the functions, $J_i$, are differentiable monotonically backward separable functions. Li showed that for problems in this class of MSOP, an equivalent multi-objective optimization problem with k-separable cost functions can be constructed. The multi-objective optimization problem can then be analytically solved, using methods relying of the differentiability of the cost function, to find the optimal input sequence for the MSOP. We do not assume, as in Li, that the cost function is differentiable or $k$-separable and our solution does not require the solution of a multi-objective optimization problem.

In related work, coherent risk measures, from~\cite{shapiro2016decomposability,shapiro2009time,ruszczynski2010risk}, result in MSOP's with non-additively separable cost functions of the form $J(\mbf u, \mbf x)= c_0(x(0),u(0))+ \rho_1 (c_1(x(1),u(1)) + \rho_2 (c_2(x(2),u(2)) + .... + \rho_T (c_T(x(T))) .... ) )$. Such MSOP's are solved recursively using a modified Bellman's Equation. Coherent risk measure functions are a special case of monotonically backward separable functions; in this case our GBE reduces to the previously proposed modified Bellman's equation.

\section{Multi-Stage Optimization Problems With Backward Separable Cost Functions} \label{sec: DP}
In this section we will introduce a class of general Multi-Stage Optimization Problems (MSOP's). We show this class contains problems that classical DP theory is able to solve; MSOP's with additively separable cost functions (Eqn.~\eqref{eqn:additively seperable}). We then propose a more general class of cost functions called monotonically backward separable functions (Eqn.~\eqref{eqn:monotonic backward seperable}) that contain the class of additively separable functions. Using this framework we are then able to derive necessary and sufficient conditions for an input sequence to solve an MSOP with monotonically backward separable cost function. Such conditions are shown to reduce to the classical conditions proposed by Bellman \cite{bellman1966dynamic} in the special case when the cost function is additively separable.



\begin{defn}
	 For a given initial condition $x_0 \in \R^n$, for every tuple of the form $\{J, f , \{X_t\}_{0 \le t \le T} , U ,T\}$, where $J: \R^{m \times T} \times \R^{n \times (T+1)} \to \R$, $f: \R^n \times \R^m \times \N \to \R^n$, $X_t \subset \R^n$, $U \subset \R^m$, and $T \in \N$, we associate a MSOP of the following form
	 \begin{align} \label{opt: DP}
	&(\mathbf{u}^*,\mathbf{x}^*) {\in} \arg \min_{\mathbf u, \mathbf x} J(\mathbf u, \mathbf x) \text{   subject to: }  \\ \nonumber
	& x(t+1)=f(x(t),u(t),t) \text{ for  } t={0},..,T-1 \\ \nonumber
	& x(0)=x_0 , \text{ } x(t) \in X_t \subset \mathbb{R}^n \text{ for  } t={0},..,T \\ \nonumber
	&u(t) \in U \subset \mathbb{R}^m \text{ for  } t={0},..,T-1\\ \nonumber
	&\mbf u=(u(0),...,u(T-1)) \text{ and } \mbf x =(x(0),...,x(T))
	\end{align}
\end{defn}
\vspace{-0.5cm}

Classical DP theory is concerned with the special case when the cost function, $J: \R^{m \times T} \times \R^{n \times (T+1)} \to \R$, has an additively separable structure defined as follows.

 \begin{defn} \label{defn: additively seperbale function}
	The function $J: U^T \times \Pi_{t=0}^T X_t \to \R$ is said to be \textbf{additively separable} if there exists functions, $c_T(x):X_T \to \mbb R$, and $c_t:X_t \times U \to \R $ for $t=0,\cdots T-1$ such that,
	
	\vspace{-1cm}
	\begin{align} \label{eqn:additively seperable}
	J(\mbf u, \mbf x)= \sum_{t=0}^{T-1} c_t(x(t),u(t)) + c_T(x(T)),
	\end{align}
	where $\mbf u=(u(0),...,u(T-1)) \text{ and } \mbf x =(x(0),...,x(T))$.
\end{defn}

We consider the class of ``monotonic backward separable'' cost functions defined next. The definition of this class of functions uses the image set of a function. Specifically, for $f:X \to Y$ we denote the image set of the function as $Image\{f\} \hspace{-0.05cm} := \hspace{-0.05cm} \{y \in Y: \hspace{-0.05cm} \text{ there exists } \hspace{-0.05cm} x \in X \hspace{-0.05cm} \text{ such that } \hspace{-0.05cm} f(x)=y\}$.

 \begin{defn} \label{defn: montonic backward seperbale function}
	The function $J: U^T \times \Pi_{t=0}^T X_t \to \R$, where $U \subset \R^m$ and $X_t \subset \R^n$ is said to be \textbf{monotonically backward separable} if there exists representation maps, $\phi_T:X_T \to \mbb R$, and $\phi_t: X_t \times U \times Image\{ \phi_{t+1} \}  \to \R $ for $t=0,\cdots T-1$ such that the following holds:
	\begin{enumerate}
		
	 	\item The function $J$ can be expressed as the composition of representation maps, $\{\phi_t\}_{t=0}^T$, ordered backwards in time. That is $J$ satisfies 
	 	
	 	\vspace{-0.6cm}
	 	{ \small \begin{align} \label{eqn:monotonic backward seperable}
 \hspace{-0.2cm} J(\mbf u, \mbf x)= \phi_0(x(0),u(0), \phi_1(x(1),u(1), \dots \phi_T(x(T)) \dots )),
		\end{align} }
		where $\mbf u=(u(0),...,u(T-1)) \text{ and } \mbf x =(x(0),...,x(T))$.
		
		\item Each representation map, $\phi_t$, is monotonic in its third argument. That is if $z,w \in Image\{\phi_{t+1}\}$ are such that $z \ge w$ \text{ then } \begin{align}\label{eqn: monotonic property}
		\phi_t(x,u,z) \ge \phi_t(x,u,w) \text{ for all } (x,u) \in X_t \times U
		\end{align}
	\end{enumerate}
Moreover if $J$ also satisfies the following properties than we say $J$ is \textbf{naturally monotonically backward separable}:
	\begin{enumerate}
		\item Each representation map, $\phi_t$, is upper semi-continuous in its third argument. That is for any $t \in \{0,..,T-1\}$, $x \in X_t$, $u \in U$ and any monotonically decreasing sequence $\{z_n\}_{n \in \N}  \subset Image\{\phi_{t+1}\}$, such that $z_{n+1} \le z_{n}$ for all $n \in \N$, then
		\begin{align} \label{property: semi cts}
		& \lim_{n\to \infty} \phi_t(x,u,z_{n}) = \phi_t(x,u, \lim_{n \to \infty} z_{n})  .
		\end{align}
		
		\item Each representation map, $\phi_t$, satisfies the following boundedness property. For any $t \in \{0,...,T-1\}$ and $(x,u,z) \in X_t \times U \times Image\{\phi_{t+1} \}$ we have $|\phi_t(x,u,z)|< \infty$ and for all $x \in X_T$ we have $|\phi_T(x)|< \infty$; that is for each $t \in \{0,...,T\}$ there exists $R>0$ such that
		\begin{align} \label{property: boundednes}
		Image \{\phi_t \} \subset \{x \in \R: |x|< R\}.
		\end{align}
	\end{enumerate}
	
\end{defn}
We show in Sec.~\ref{sec: POP} that monotonically backward separable functions share a deep connection with Bellman's Principle of Optimality (Defn.~\ref{defn: principle of optimality deterministic}). However, we also consider naturally monotonically backward separable functions as the added semi-continuity and boundedness properties are used in the derivation of necessary and sufficient conditions for an input sequence to solve an MSOP with naturally monotonically backward separable cost function (Theorem~\ref{thm: GBE}).

We next show the class of MSOP's with monotonically backward separable cost functions includes the class of MSOP's with additively separable cost functions as a special case.

\begin{lem} \label{lem: add sep is back sep}
	 Suppose $J: U^T \times \Pi_{t=0}^T X_t \to \R$ is an additively separable function (Defn.~\ref{defn: additively seperbale function}), with associated cost functions $\{c_t\}_{t=0}^T$. Then $J$ is a monotonically backward separable function (Defn.~\ref{defn: montonic backward seperbale function}). Moreover, if the functions $\{c_t\}_{t=0}^T$ are bounded over $X_t \times U$ then $J$ is naturally monotonically backward separable function.
\end{lem}
\begin{pf}
	Given an additively separable function, $J$, we know there exists functions $\{c_t\}_{0 \le t \le T}$ such that Eqn.~\eqref{eqn:additively seperable} holds. To prove $J$ is monotonically backward separable we construct representation maps $\{\phi_t\}_{t=0}^T$ such that Eqns.~\eqref{eqn:monotonic backward seperable} and \eqref{eqn: monotonic property} hold. We define these representation maps as follows:
	\begin{align} \label{eqn: addative functions are backward sep}
	& \phi_i(x,u,z)=c_i(x,u)+z \quad \text{for } i=0,\cdots,T-1 \\ \nonumber
	& \phi_T(x,w)=c_T(x).
	\end{align}
	Now, $\frac{\partial \phi_t(x,y,z)}{\partial z}=1>0$ for all $t \in \{0,....,T-1\}$, $x \in X_t$ and $u \in U$, implying the monotonicity property in Eqn.~\eqref{eqn: monotonic property}. 
	
	Now assuming the functions $\{c_t\}_{t=0}^T$ are bounded over $X_t \times U$ it follows trivially that the representation maps $\{\phi_t\}_{t=0}^T$, given in Eqn.~\eqref{eqn: addative functions are backward sep},  satisfy the semi-continuity and boundedness properties given in Eqns.~\eqref{property: semi cts} and \eqref{property: boundednes}. Thus $J$ is naturally monotonically backward separable function. $\blacksquare$
	\end{pf}
Further examples of monotonically backward separable functions, including instances where the representation maps are non-differentiable, are given in Section \ref{subsec: Examples}.

\subsection{Exchanging The Order Of Composition And Infimum For Monotonically Backward Separable Functions}
As we will show in Lemma \ref{lem: interchanging limits}, monotonically backward separable functions have the special property that the order of an infimum and composition of representation maps can be interchanged. To show this we must use the monotonic convergence theorem.
\begin{thm}[Monotone Convergence Theorem] \label{thm: monotone convergence theorem}
	Suppose $\{z_{n}\}_{n \in \N} \subset \R$ is a bounded sequence that is monotonically decreasing, $z_{n+1} \le z_n$ for all $n \in \N$. Then $\lim_{n \to \infty}z_n = \inf_{n \in \N} z_n$.
\end{thm}
Before proving in Lemma \ref{lem: interchanging limits} we introduce notation for the set of feasible controls. Given a tuple $\{J, f , \{X_t\}_{0 \le t \le T} , U ,T\}$ for $x \in X_t$ and $s \in [0,T-1]$ we denote
\begin{align*}
\Gamma_{x,s}:=\{u \in U: f(x,u,s) \in X_{s+1}\}.
\end{align*}
Moreover we say
\begin{equation} \label{notation: feasible input sequence}
(u(s),...,u(T-1)) \in \Gamma_{x_0,[s,T-1]}
\end{equation}
if $u(t) \in \Gamma_{x(t),t}$ for all $t \in \{s,...,T-1\}$, where $x(s)=x_0$ and $x(k+1)= f(x(k),u(k),k)$ for $k \in \{s,...,T-1\}$.
\begin{lem} \label{lem: interchanging limits}
{ \normalsize Consider an MSOP of Form~\eqref{opt: DP} associated with $\{J, f , \{X_t\}_{0 \le t \le T} , U ,T\}$. Suppose $J: U^T \times \Pi_{t=0}^T X_t \to \R$ is a naturally monotonically backward separable function (Defn.~\ref{defn: montonic backward seperbale function}) with representation maps $\{\phi_t\}_{t=0}^T$ and $\Gamma_{x,t} \ne \emptyset$ for all $(x,t) \in X_t \times \{0,...,T-1\}$. Then for $k \in \{0,...,T-1\}$ and any $x \in X_k$ we have }

\vspace{-0.6cm}
{\small \begin{align} \nonumber
	&  \inf_{u(k) \in \Gamma_{x,k}} \bigg\{  \phi_k\bigg(x(k),u(k), \inf_{(u(k+1),...,u(T-1)) \in \Gamma_{x(k+1),[k+1,T-1]}} \bigg\{ \phi_{k+1}( \\ \nonumber
& x(k+1),u(k+1), \phi_{k+2}(x(k+2),u(k+2),... \phi_T(x(T))...))  \bigg\} \bigg) \bigg\}\\ \label{eqn: interchange inf}
& = \inf_{(u(k),...,u(T-1)) \in \Gamma_{x,[k,T-1]}} \bigg\{ \\ \nonumber
& \qquad \qquad \phi_k(x(k),u(k), \phi_{k+1}(x(k+1),u(k+1),...  \phi_T(x(T))...)) \bigg\},
\end{align} }
\noindent where $x(t+1)=f(x(t),u(t),t)$ for $t \in \{k,...,T-1\}$ and $x(k)=x$.
\end{lem}
\begin{pf} To show Eqn.~\eqref{eqn: interchange inf} we will split the proof into two parts. In Part 1 we will show the left hand side of Eqn.~\eqref{eqn: interchange inf} is less than or equal to the right hand side of Eqn.~\eqref{eqn: interchange inf}. In Part 2 we will show the right hand side of Eqn.~\eqref{eqn: interchange inf} is less than or equal to the left hand side of Eqn.~\eqref{eqn: interchange inf}.
	
\textbf{Part 1 of proof:} By the definition of the infimum it follows for all $y \in X_{k+1}$ that

\vspace{-0.6cm}
{ \small \begin{align} \nonumber
& \inf_{(u(k+1),...,u(T-1)) \in \Gamma_{y,[k+1,T-1]}} \phi_{k+1}({x}(k+1),u(k+1),...  \phi_T(x(T))...) \\ \label{inf}
& \le \phi_{k+1}(\tilde{x}(k+1),\tilde{u}(k+1),...  \phi_T(\tilde{x}(T))...),
\end{align} } for any $(\tilde{u}(k+1),...,\tilde{u}(T-1)) \in \Gamma_{x(k+1),[k+1,T-1]}$, where $\tilde{x}(t+1)=f(\tilde{x}(t),\tilde{u}(t),t)$, ${x}(t+1)=f({x}(t),{u}(t),t)$ for $t \in \{k+1,...T-1\}$, and $x(k+1)=\tilde{x}(k+1)=y$.

Since $\phi_k$ is monotonic in its third argument (Eqn.~\eqref{eqn: monotonic property}) it follows from Eqn.~\eqref{inf} that for any $(x,u) \in X_k \times \Gamma_{x,k}$ that
\begin{align} \label{2}
& \phi_k(x(k),u(k), \inf_{(u(k+1),...,u(T-1)) \in \Gamma_{x(k+1),[k+1,T-1]}}\{\\ \nonumber
&  \qquad \qquad \phi_{k+1}({x}(k+1),u(k+1),...  \phi_T(x(T))...)  \}\\ \nonumber
& \le \phi_k(x(k),u(k),\phi_{k+1}(\tilde{x}(k+1),\tilde{u}(k+1),...  \phi_T(\tilde{x}(T))...) ),
\end{align}
for any $(\tilde{u}(k+1),...,\tilde{u}(T-1)) \in \Gamma_{x(k+1),[k+1,T-1]}$, where $\tilde{x}(t+1)=f(\tilde{x}(t),\tilde{u}(t),t)$, ${x}(t+1)=f({x}(t),{u}(t),t)$ for $t \in \{k,...T-1\}$, $x(k)=\tilde{x}(k)=x$, and $u(k)=u$.

Now, since Eqn.~\eqref{2} holds for any $u \in  \Gamma_{x,k}$ and $(\tilde{u}(k+1),...,\tilde{u}(T-1)) \in \Gamma_{x(k+1),[k+1,T-1]}$ we are able to take the infimum over these in Eqn.~\eqref{2}, deducing the left hand side of Eqn.~\eqref{eqn: interchange inf} is less or than or equal to its right hand side.

\textbf{Part 2 of proof:} Let us fix $(x,u) \in X_k \times \Gamma_{x(k),k}$. Since $\Gamma_{x,t} \ne \emptyset$ for all $(x,t) \in X_t \times \{0,...,T-1\}$ it follows from the definition of the infimum for all $n \in \N$ there exists $({u}_n(k+1),...,{u}_n(T-1)) \in \Gamma_{x(k+1),[k+1,T-1]}$ such that

\vspace{-0.6cm} 
{ \small\begin{align} \nonumber
& \inf_{(u(k+1),...,u(T-1)) \in \Gamma_{x(k+1),[k+1,T-1]}}   \hspace{-0.5cm} \phi_{k+1}({x}(k+1),u(k+1),...  \phi_T(x(T))...)\\ \label{111}
& \le \phi_{k+1}({x}_n(k+1),u_n(k+1),...  \phi_T(x_n(T))...) \\ \nonumber
&  \le \hspace{-0.1cm} \inf_{(u(k+1),...,u(T-1)) \in \Gamma_{x(k+1),[k+1,T-1]}} \hspace{-1.4cm} \phi_{k+1}({x}(k+1),u(k+1),...  \phi_T(x(T))...) + \frac{1}{n},
\end{align} }
{\noindent where $x_n(t+1)=f(x_n(t),u_n(t),t)$ for $t \in \{k+1,...,T-1\}$, and $x_n(k+1)=x(k+1)=f(x,u,k)$. \noindent }

Let $a_n:=\phi_{k+1}({x}_n(k+1),u_n(k+1),...  \phi_T(x_n(T))...)$. It follows from Eqn.~\eqref{111} that

\vspace{-0.6cm}
{\small \begin{align*}
&\lim_{n \to \infty} a_n = \hspace{-0.4cm} \inf_{(u(k+1),...,u(T-1)) \in \Gamma_{x(k+1),[k+1,T-1]}} \hspace{-1.4cm} \phi_{k+1}({x}(k+1),u(k+1),...  \phi_T(x(T))...).\\
&a_n \ge \inf_{(u(k+1),...,u(T-1)) \in \Gamma_{x(k+1),[k+1,T-1]}} \hspace{-1.4cm} \phi_{k+1}({x}(k+1),u(k+1),...  \phi_T(x(T))...),\\
& \hspace{6cm} \text{for all } n \in \N.
\end{align*} }
Since $\{a_n\}_{n \in \N}$ converges to some limit from above there exists a monotonically decreasing subsequence $\{b_n\}_{n \in \N} \subseteq \{a_n\}_{n \in \N}$ such that $b_{n+1} \le b_{n}$ for $n \in \N$. Using $\{b_n\}_{n \in \N}$ we now define
\[z_n:= \phi_k(x,u,b_n).\]
Since $\phi_k$ is monotonic in its third argument (Eqn.~\eqref{eqn: monotonic property}) and $b_{n+1} \le b_{n}$ it follows $z_{n+1}=\phi_k(x,u,b_{n+1}) \le \phi_k(x,u,b_{n}) \le z_{n}$. Hence $\{z_n\}_{n \in \N}$ is a monotonically decreasing sequence. Moreover, since $\phi_k$ has the property that it is a bounded over $X_k \times U \times  Image\{\phi_{k+1}\}$ (Eqn.~\eqref{property: boundednes}) it follows that $\{z_n\}_{n \in \N}$ is a bounded sequence. Now by the monotone convergence theorem (Theorem \ref{thm: monotone convergence theorem}) we have that $\inf_{n \in \N} z_n = \lim_{n \to \infty} z_n$.

It now follows since $\phi_k$ is upper semi-continuous (Eqn.~\eqref{property: semi cts}) in its third argument that
\begin{align} \label{4}
&  \inf_{(u(k+1),...,u(T-1)) \in \Gamma_{x(k+1),[k,T-1]}} \{ \\ \nonumber
& \qquad \qquad \phi_k(x,u, \phi_{k+1}(x(k+1),u(k+1),...  \phi_T(x(T))...)) \}\\ \nonumber
& \le \inf_{n \in \N} \phi_k(x,u, \phi_{k+1}(x_n(k+1),u_n(k+1),...  \phi_T(x_n(T))...))\\ \nonumber
&\le \inf_{n \in \N} z_n = \lim_{n \to \infty} z_n\\ \nonumber
& =\lim_{n \to \infty} \phi_k(x,u, b_n) = \phi_k(x,u, \lim_{n \to \infty} b_n) = \phi_k(x,u, \lim_{n \to \infty} a_n)\\ \nonumber
&= \phi_k(x,u, \inf_{(u(k+1),...,u(T-1)) \in \Gamma_{x(k+1),[k+1,T-1]}}\{\\ \nonumber
&  \qquad \qquad \phi_{k+1}({x}(k+1),u(k+1),...  \phi_T(x(T))...)  \}.
\end{align}
Since Eqn.~\eqref{4} holds for any arbitrarily selected $(x,u) \in X_k \times \Gamma_{x,k}$ we are able to take the infimum with respect to $u \in \Gamma_{x,k}$, showing the right hand side of Eqn.~\eqref{eqn: interchange inf} is less than or equal to its left hand side. 

In Part 1 of the proof we have shown that the left hand side of Eqn.~\eqref{eqn: interchange inf} is less than or equal to the right hand side of Eqn.~\eqref{eqn: interchange inf}. In Part 2 of the proof we have shown that the right hand side of Eqn.~\eqref{eqn: interchange inf} is less than or equal to the left hand side of Eqn.~\eqref{eqn: interchange inf}. Putting these two parts together we deduce the left hand side must equal the right hand side, therefore completing the proof and showing Eqn.~\eqref{eqn: interchange inf} holds. $\blacksquare$ \end{pf}


\subsection{Main Result: A Generalization Of Bellman's Equation}\label{subsec:main}

When $J$ is additively separable, the MSOP, given in Eqn.~\eqref{opt: DP}, associated with the tuple $\{J, f , \{X_t\}_{0 \le t \le T} , U ,T\}$, can be solved recursively using Bellman's Equation~\cite{bellman1966dynamic}. In this section we show that a similar approach can be used to solve MSOP's with naturally monotonically backward separable cost functions. 

We next define conditions under which a function, $V$, is said to be a \textit{value function} for an associated MSOP.
\begin{defn} \label{defn: Value functions}
	Consider a monotonically backward separable function  $J: \R^{m \times T} \times \R^{n \times (T+1)} \to \R$ with representation functions $\{\phi_t\}_{0 \le t \le T}$, $f: \R^n \times \R^m \times \N \to \R^n$, $X \subset \R^n$, $U \subset \R^m$, and $T \in \N$. We say the function $V: \R^n \times [0,T] \to \R$ is a \textbf{value function} of the MSOP associated with the tuple $\{J, f , \{X_t\}_{0 \le t \le T} , U ,T\}$ if for all $x \in X_T$
	\begin{align} \label{fun: value 1}
	V(x,T)= \phi_T(x),
	\end{align}
and for all $x \in X_t$  and $t \in [0,T-1]$
\begin{align} \label{fun: value 2}
&V(x,t)= \inf_{u(t) \in \Gamma_{x,t},...., u(T-1) \in \Gamma_{x(T-1),T-1}} \bigg\{ \\ \nonumber
& \phi_t(x(t),u(t), \phi_{t+1}(x(t+1),u(t+1),... \phi_T(x(T))...)) \bigg\},
\end{align}
where $x(t)=x$ and $x(k+1)= f(x(k),u(k),k)$ for $k \in \{t,...,T-1\}$.
\end{defn}
We note that the value function has the special property that $V(x_0,0)=J^*$, where $J^*$ is the minimum value of the cost function of the MSOP \eqref{opt: DP}. In the special case when $J$ is an additively separable function the value function defined in this way reduces to the optimal cost-to-go function.

\begin{prop}[Generalized Bellman's Equation (GBE)] \label{prop: gen bellman}
	 Consider an MSOP of Form~\eqref{opt: DP} associated with \linebreak $\{J, f , \{X_t\}_{0 \le t \le T} , U ,T\}$. Suppose $J: U^T \times \Pi_{t=0}^T X_t \to \R$ is a naturally monotonically backward separable function (Defn.~\ref{defn: montonic backward seperbale function}) with representation maps $\{\phi_t\}_{t=0}^T$ and $\Gamma_{x,t} \ne \emptyset$ for all $(x,t) \in X_t \times \{0,...,T-1\}$. Then if $F: \R^n \times [0,T] \to \R$ satisfies
	\begin{align} \label{eqn: generalize Bellman}
	F(x,T) & = \phi_T(x) \text{    for all } x \in X_T, \quad \text{and} \\ \nonumber
	 F(x,t) & = \inf_{u \in \Gamma_{x,t}} \bigg\{ \phi_t(x,u,F(f(x,u,t),t+1)) \bigg\} \\ \nonumber
	 & \qquad \quad \text{ for all } x \in X_t, t \in \{0,..,T-1\},
	\end{align}
 then $F$ is a value function (Defn.~\ref{defn: Value functions}) of the MSOP associated with $\{J, f , \{X_t\}_{0 \le t \le T} , U ,T\}$.
\end{prop}

\begin{pf}
	Suppose $F$ satisfies Eqn.~\eqref{eqn: generalize Bellman}. To show $F$ is a value function of the MSOP associated with the tuple $\{J, f , \{X_t\}_{0 \le t \le T} , U ,T\}$ we must show it satisfies Eqns.~\eqref{fun: value 1} and \eqref{fun: value 2}. We prove this using backward induction in the time variable of $F$. Clearly $F$ satisfies Eqn.~\eqref{fun: value 1} for $k=T$. Now, for our induction hypothesis, let us assume for some $k \in \{0,...,T-1\}$ that $F$ satisfies Eqn.~\eqref{fun: value 2} at time-stage $k+1$ for all $x \in X_{k+1}$. We will now show that the induction hypothesis implies $F$ must also satisfy Eqn.~\eqref{fun: value 2} at time-stage $k$ for all $x \in X_{k}$. Letting $x \in X_{k}$ we have 	
	\vspace{-0.6cm}
	{\small
		\begin{align*}
	& F(x,k)  = \inf_{u \in \Gamma_{x,k}} \bigg\{ \phi_k(x,u,F(f(x,u,k),k+1)) \bigg\}\\
	& = \inf_{u \in \Gamma_{x,k}} \bigg\{  \phi_k\bigg(x,u, \inf_{u(k+1) \in \Gamma_{x(k+1),k+1},...., u(T-1) \in \Gamma_{x(T-1),T-1}} \bigg\{ \phi_{k+1}( \\ \nonumber
	& x(k+1),u(k+1), \phi_{k+2}(x(k+2),u(k+2),... \phi_T(x(T))...))  \bigg\} \bigg) \bigg\}\\
	& = \inf_{u(k) \in \Gamma_{x,k},...., u(T-1) \in \Gamma_{x(T-1),T-1}} \bigg\{ \\ \nonumber
	& \qquad \qquad \phi_k(x(k),u(k), \phi_{k+1}(x(k+1),u(k+1),...  \phi_T(x(T))...)) \bigg\},
	\end{align*} } %
\vspace{-0.75cm}

{where $x(k)=x$ and $x(t+1)= f(x(t),u(t),t)$ for $t \in \{k,...,T-1\}$. The first equality follows as $F$ satisfies Eqn.~\eqref{eqn: generalize Bellman}; the second equality follows from the induction hypothesis; the third equality follows by Lemma~\ref{lem: interchanging limits}. \noindent }

Therefore, by backward induction, we conclude $F$ satisfies Eqns.~\eqref{fun: value 1} and \eqref{fun: value 2} and hence is a value function for the MSOP associated with the tuple $\{J, f , \{X_t\}_{0 \le t \le T} , U ,T\}$. $\blacksquare$ \end{pf}

{We next propose sufficient conditions showing an input sequence is optimal if it recursively minimizes the right hand side of the GBE~\eqref{eqn: generalize Bellman}. Later in Theorem~\ref{thm: GBE} we propose necessary and sufficient conditions involving the GBE~\eqref{eqn: generalize Bellman}.}

\begin{prop}[Sufficient conditions]
	 \label{prop: suff GBE}
	Consider an MSOP of Form~\eqref{opt: DP} associated with $\{J, f , \{X_t\}_{0 \le t \le T} , U ,T\}$. Suppose $J: U^T \times \Pi_{t=0}^T X_t \to \R$ is a naturally monotonically backward separable function (Defn.~\ref{defn: montonic backward seperbale function}) with representation maps $\{\phi_t\}_{t=0}^T$, $\Gamma_{x,t} \ne \emptyset$ for all $(x,t) \in X_t \times \{0,...,T-1\}$, $V: \R^n \times [0,T] \to \R$ satisfies the GBE~\eqref{eqn: generalize Bellman}, and the state sequence $\mathbf x^*=(x^*(0),...,x^*(T))$ and input sequence $\mathbf u^*=(u^*(0),...,u^*(T-1))$ satisfy
		\begin{align} \nonumber 
		& u^*(k) \in   \arg \inf_{ u \in \Gamma_{x^*(k),k}} \bigg\{ \phi_t(x^*(k),u,V(f(x^*(k),u,k),k+1)) \bigg\}\\ \label{fun: optimal policy}
		& \hspace {1cm} \text{ for } k\in \{0,...,T-1\}.\\  \nonumber
				 & x^*(0)=x_0, \quad x^*(k+1)= f(x^*(k),u^*(k),k)\\ \label{iff 1}
		& \hspace {1cm} \text{ for } k\in \{0,...,T-1\}.
		\end{align}
Then $(\mbf u^*, \mbf x^*)$ solve the MSOP given in Eqn.~\eqref{opt: DP}, associated with the tuple $\{J, f , \{X_t\}_{0 \le t \le T} , U ,T\}$.
\end{prop}

\begin{pf}
%
Suppose $(\mathbf u^*,\mathbf x^*)$ satisfy Eqns.~\eqref{fun: optimal policy} and \eqref{iff 1}. It follows the pair ($\mathbf u^*$,$\mathbf x^*$) is a feasible solution for MSOP given in Eqn.~\eqref{opt: DP} since Eqn.~\eqref{fun: optimal policy} implies $u^*(k) \in \Gamma_{x^*(k),k}$, thus $u^*(k) \in U$ and, using Eqn.~\eqref{iff 1}, $x^*(k+1)=f(x^*(k),u^*(k),k) \in X_{k+1}$ for all $k \in \{0,...,T-1\}$.

By Eqn.~\eqref{fun: optimal policy} it follows for all $k \in \{0,...,T-1\}$ that
\begin{align} \label{u is optimal}
\inf_{ u \in \Gamma_{x^*(k),k}} \bigg\{ & \phi_k(x^*(k),u,V(f(x^*(k),u,k),k+1)) \bigg\} \\ \nonumber
& =\phi_k(x^*(k),u^*(k),V(f(x^*(k),u^*(k),k),k+1)).
\end{align}
We will now show Eqn.~\eqref{u is optimal} implies $(\mathbf u^*,\mathbf x^*)$ solve the MSOP. 
\vspace{-0.75cm}
{\small \begin{align*}
& \inf_{\mathbf u \in \Gamma_{x_0,[0,T-1]}} J(\mathbf u, \mathbf x) =V(x_0,0)\\
& = \inf_{u \in \Gamma_{x^*(0),0}} \bigg\{ \phi_0(x^*(0),u,V(f(x^*(0),u,0),1)) \bigg\} \\
& = \phi_0(x^*(0),u^*(0),V(x^*(1),1))\\
&= \phi_0\bigg(x^*(0),u^*(0), \inf_{u \in \Gamma_{x^*(1),1}} \bigg\{ \phi_1(x^*(1),u,V(f(x^*(1),u,1),2))  \bigg\} \bigg)\\
& \hspace{3cm} \vdots\\
&= \phi_0(x^*(0),u^*(0),..., \phi_k(x^*(k),u^*(k),\\
& \hspace{2cm} \phi_{k+1}(x^*(k+1),u^*(k+1),.... \phi_T(x^*(T)))...)...)\\
&= J(\mathbf u^*, \mathbf x^*),
\end{align*}}
\vspace{-1cm}

{where the first equality follows as it was shown in Prop.~\ref{prop: gen bellman} that $V$ is a value function of the MSOP, the second equality follows since $V$ satisfies the GBE~\eqref{eqn: generalize Bellman} and using $x^*(0)=x_0$, the third equality follows by Eqn.~\eqref{u is optimal}, the fourth inequality follows again using the GBE, and the fifth inequality follows by recursively using the GBE together with Eqn.~\eqref{u is optimal}. Thus if $( \mathbf u^*, \mathbf x^*)$ satisfy Eqns.~\eqref{iff 1} and \eqref{fun: optimal policy} then $ (\mathbf u^*, \mathbf x^*)$ solve the MSOP given in Eqn.~\eqref{opt: DP}. $\blacksquare$  \noindent}
\end{pf}
Consider an MSOP associated with $\{J, f , \{X_t\}_{0 \le t \le T} , U ,T\}$, where $J$ is naturally monotonically backward separable (Defn.~\ref{defn: montonic backward seperbale function}). As we will show next, if the representation maps $\{\phi_t\}_{t=0}^T$, associated with $J$ are strictly monotonic (Eqn.~\eqref{strict monotonic}) then Eqns.~\eqref{fun: optimal policy} and \eqref{iff 1} of Prop.~\ref{prop: suff GBE} become sufficient and necessary for optimality. In Sec.~\ref{subsec: Examples} we will give several examples of naturally monotonically backward functions with associated strictly monotonic representation maps.
\begin{thm}[Necessary and sufficient conditions] \label{thm: GBE}
Consider an MSOP of Form~\eqref{opt: DP} associated with \linebreak $\{J, f , \{X_t\}_{0 \le t \le T} , U ,T\}$. Suppose $J: U^T \times \Pi_{t=0}^T X_t \to \R$ is a naturally monotonically backward separable function (Defn.~\ref{defn: montonic backward seperbale function}) with representation maps $\{\phi_t\}_{t=0}^T$, and $\Gamma_{x,t} \ne \emptyset$ for all $(x,t) \in X_t \times \{0,...,T-1\}$. Furthermore, suppose the representation maps are strictly monotonic in their third argument. That is if $z,w \in Image\{\phi_{t+1}\}$ are such that $z > w$ then  \begin{align} \label{strict monotonic}
\phi_t(x,u,z) > \phi_t(x,u,w) \text{ for all } (x,u) \in X_t \times U.
\end{align}
Then $(\mathbf u^*, \mathbf x^*)$ solve the MSOP if and only if $ (\mathbf u^*, \mathbf x^*)$ satisfy Eqns.~\eqref{fun: optimal policy} and \eqref{iff 1}.
\end{thm}
\begin{pf}
 If $(\mathbf u^*,\mathbf x^*)$ satisfy Eqns.~\eqref{fun: optimal policy} and \eqref{iff 1} then Prop.~\ref{prop: suff GBE} shows $ (\mathbf u^*, \mathbf x^*)$ solve the MSOP associated with $\{J, f , \{X_t\}_{0 \le t \le T} , U ,T\}$.
 
    Now assume the representation maps $\{\phi_t\}_{t=0}^T$ are strictly monotonic in their third argument (Eqn.~\eqref{strict monotonic}) and $(\mathbf u^*,\mathbf x^*)$ solve the MSOP associated with the tuple $\{J, f , \{X_t\}_{0 \le t \le T} , U ,T\}$. As we have assumed $(\mathbf u^*,\mathbf x^*)$ is a solution then it follows $(\mathbf u^*,\mathbf x^*)$ is feasible and thus Eqn.~\eqref{iff 1} is trivially satisfied. To prove Eqn.~\eqref{fun: optimal policy} is also satisfied let us suppose for contradiction the negation of Eqn.~\eqref{fun: optimal policy}, that there exists $k \in \{0,...,T-1\}$ such that
\begin{align*}
u^*(k) \notin   \arg \inf_{ u \in \Gamma_{x^*(k),k}} \bigg\{ \phi_t(x^*(k),u,V(f(x^*(k),u,k),k+1)) \bigg\},
\end{align*}
where $V: \R^n \times [0,T] \to \R$ satisfies the GBE~\eqref{eqn: generalize Bellman}, and hence it follows
\begin{align} \label{contradiction: assum}
\inf_{ u \in \Gamma_{x,t}} \bigg\{ & \phi_k(x^*(k),u,V(f(x^*(k),u,k),k+1)) \bigg\} \\ \nonumber
& < \phi_k(x^*(k),u^*(k),V(f(x^*(k),u^*(k),k),k+1)).
\end{align}
Using Eqn.~\eqref{contradiction: assum} it follows,
\vspace{-0.75cm} { \small \begin{align*}
	& J(\mathbf u^*, \mathbf x^*)= \inf_{\mathbf u \in \Gamma_{x_0,[0,T-1]}} J((u(0),..,u(T-1)),(x(0),...,x(T))) \\
	& \le \inf_{\mathbf w \in \Gamma_{x^*(k),[k,T-1]} } J((u^*(0),..,u^*(k-1),w(k),..,w(T-1)), \\ \nonumber
	& \hspace{3.5cm} (x^*(0),...,x^*(k),z(k+1),...,z(T))) \\ \nonumber
	& = \phi_0 \bigg( x^*(0),u^*(0),..., \inf_{w(k) \in \Gamma_{x^*(k),k}} \bigg\{ \phi_{k}(x^*(k),w(k), \\
	&   \inf_{\mathbf w \in \Gamma_{f(x^*(k),w(k),k),[k+1,T-1]}}\phi_{k+1}(z(k+1),w(k+1),.... \phi_T(z(T))...) \bigg\}... \bigg)  \\
	&=\phi_0 \bigg(x^*(0),u^*(0),...,\\
	& \quad \inf_{w(k) \in \Gamma_{x^*(k),k}} \bigg\{ \phi_{k}(x^*(k),w(k),V(f(x^*(k),w(k),k),k+1)) \bigg\},..,\bigg) \\
	& < \phi_0(x^*(0),u^*(0),...,\\ 
	& \hspace{1cm} \phi_k(x^*(k),u^*(k),V(f(x^*(k),u^*(k),k),k+1)),..,) \\
	&= \phi_0\bigg(x^*(0),u^*(0),..., \phi_k \bigg(x^*(k),u^*(k), \inf_{\mathbf w \in \Gamma_{f(x^*(k),w(k),k),[k+1,T-1]}} \\
	&  \hspace{1.5cm} \bigg\{ \phi_{k+1}(z(k+1),w(k+1),... \phi_T(z(T)))...) \bigg\} \bigg)...\bigg) \\
	& \le \phi_0(x^*(0),u^*(0),..., \phi_k(x^*(k),u^*(k),\\
	& \qquad \phi_{k+1}(x^*(k+1),u^*(k+1),.... \phi_T(x^*(T)))...)...)\\
	&= J(\mathbf u^*, \mathbf x^*),
	\end{align*} } \normalsize
\vspace{-1cm}

{where the first equality follows as the pair $(\mathbf u^*\mathbf x^*)$ is assumed to solve the MSOP. The first inequality follows by taking the infimum only over the input and state sequences from time stage $k+1$ onwards and fixing the first $k$ input and state sequences as $(u^*(0),..,u^*(k-1))$ and $(x^*(0),...,x^*(k))$ (which are known to be feasible as the pair $(\mathbf u^*,\mathbf x^*)$ is assumed to solve the MSOP). The second equality follows by Lemma~\ref{lem: interchanging limits}. The third equality follows by Prop.~\ref{prop: gen bellman} that shows $V$ is the value function. The second inequality follows from Eqn.~\eqref{contradiction: assum} and using the assumed strict monotonic property of the representation maps (Eqn.~\eqref{strict monotonic}). The fourth equality follows using Prop.~\ref{prop: gen bellman}, that shows $V$ is the value function. The third inequality follows by fixing the decision variables of the infimum to $(u^*(k),...,u^*(T-1))$ and $(x^*(k+1),...,x^*(T))$ (which are known to be feasible as the pair $(\mathbf u^*,\mathbf x^*)$ is assumed to solve the MSOP) and using monotonic property of the representation maps (Eqn.~\eqref{eqn: monotonic property}). \noindent }

We therefore get a contradiction, that $J(\mathbf u^*, \mathbf x^*)<J(\mathbf u^*, \mathbf x^*)$; showing if $(\mathbf u^*,\mathbf x^*)$ solve the MSOP then Eqns.~\eqref{iff 1} and \eqref{fun: optimal policy} must hold. $\blacksquare$ \end{pf}


In the next corollary we show that when the cost function, $J$, is additively separable, the GBE~\eqref{eqn: generalize Bellman} reduces to Bellman's Equation~\eqref{eqn: Bellman Eqn}; thus showing Bellman's Equation is an implication of the GBE. Therefore we have generalized the necessary and sufficient conditions for optimality encapsulated in Bellman's Equation to the GBE. The GBE provides optimality conditions for a larger class of MSOP's with monotonically backward separable cost functions; that no longer need be additively separable.



\begin{cor}[Bellman's Equation] \label{cor: GBE implies BE}
	Consider an MSOP of Form~\eqref{opt: DP} associated with $\{J, f , \{X_t\}_{0 \le t \le T} , U ,T\}$. Suppose $J: U^T \times \Pi_{t=0}^T X_t \to \R$ is an additively separable function (Defn.~\ref{defn: additively seperbale function}), with associated cost functions $\{c_t\}_{t=0}^T$ that are bounded over $X_t \times U$. Then if $F: \R^n \times [0,T] \to \R$ satisfies
\begin{align} \label{eqn: Bellman Eqn}
F(x,T) & = c_T(x) \quad \text{ for all } x \in X_T, \\ \nonumber
F(x,t) & = \inf_{u \in \Gamma_{x,t}} \bigg\{ c_t(x,u)+F(f(x,u,t),t+1) \bigg\} \\ \nonumber
& \qquad \qquad \text{ for all } x \in X_t, \;t \in \{0,..,T-1\},
\end{align}
then $F$ is a value function for the MSOP associated with the tuple $\{J, f , \{X_t\}_{0 \le t \le T} , U ,T\}$.

Moreover, if $\Gamma_{x,t} \ne \emptyset$ for all $(x,t) \in X_t \times \{0,...,T\}$ then $\mathbf x^*=(x^*(0),...,x^*(T))$ and $\mathbf u^*=(u^*(0),...,u^*(T-1))$ solve the MSOP if and only if the following is satisfied

\vspace{-0.5cm}
{\small \begin{align} \label{22}
&u^*(k)  \in  \arg  \inf_{ u \in \Gamma_{{x}^*(k),k}} \hspace{-0.15cm} \{c_k(x^*(k),u) +F(f(x^*(k),u,k),k+1) \},\\ \label{23}
& x^*(0)=x_0, \quad x^*(k+1)= f(x^*(k),u^*(k),k)\\ \nonumber
& \hspace {1cm} \text{ for } k\in \{0,...,T-1\}.
\end{align} } \normalsize
\end{cor}
\begin{pf}
	By Lemma~\ref{lem: add sep is back sep} it follows $J$ is naturally monotonically backward separable and can be written in Form~\eqref{eqn:monotonic backward seperable} using the representation maps given in Eqn.~\eqref{eqn: addative functions are backward sep}. Substituting the representation maps in Eqn~\eqref{eqn: addative functions are backward sep} into the GBE~\eqref{eqn: generalize Bellman}, we obtain Bellman's Equation~\eqref{eqn: Bellman Eqn}. Prop.~\ref{prop: gen bellman} then shows $F$ is a value function for the MSOP, associated with the tuple $\{J, f , \{X_t\}_{0 \le t \le T} , U ,T\}$. 
	
	Moreover as the representation maps in Eqn.~\eqref{eqn: addative functions are backward sep} are clearly strictly monotonic in their third argument (Eqn.~\eqref{strict monotonic}) it follows by Theorem~\ref{thm: GBE} that $(\mbf x^*, \mbf u^*)$ solve the MSOP if and only if $(\mbf x^*, \mbf u^*)$ satisfy Eqns.~\eqref{22} and \eqref{23}. $\blacksquare$
	\end{pf}

\normalsize
\subsection{Examples: Backward Separable Functions} \label{subsec: Examples}
In Subsection~\ref{subsec:main}, we have shown that MSOP's with cost functions that are naturally monotonically backward separable (Defn.~\ref{defn: montonic backward seperbale function}) can be solved efficiently using the GBE~\eqref{eqn: generalize Bellman}. We now give examples of non-additively separable, yet monotonically backward separable functions, which may be of significant interest. This is not a complete list of all monotonically backward separable functions. Currently little is known about size and structure of the set of all monotonically backward separable functions.

%

The first function we consider is the point-wise maximum function. This function occurs in MSOP's when demand charges are present~\cite{jones2018extensions} and in maximal invariant set estimation~\cite{xue2018robust}. 
\begin{ex} [Point wise maximum function] \label{ex: sup forward sep} Suppose $J: U^T \times \Pi_{t=0}^T X_t \to \R$ is of the form
	\[
	J(\mbf u, \mbf x)=\max \left\{\max_{0 \le k \le T-1}\{c_k(x(k),u(k))\}, c_T(x(T)) \right\},
	\]
	where $\mbf u=(u(0),...,u(T-1))$, $\mbf x =(x(0),...,x(T))$, $U \subseteq \R^m$, $X_t \subseteq \R^n$, $c_k: X_k \times U \to \R$ and $c_T:X_T \to \R$. Then $J$ is a monotonically backward separable function. Moreover, if $\{c_t\}_{t=0}^T$ are bounded functions, then $J$ is naturally monotonically backward separable.
\end{ex}
\begin{pf} We can write $J$ in Form~\eqref{eqn:monotonic backward seperable} using the representation functions
	{\begin{align} \label{rep: pointwise max}
		\phi_T(x) & = c_T(x),\\ \nonumber
		 \phi_{i}(x,u,z) & =\max\{c_{i}(x,u), z\}  
		\text{ for all } i \in \{0,..,T-1\}.
		\end{align} }
	The monotonicity property in Eqn.~\eqref{eqn: monotonic property} follows since if $y \ge z$ then for all $ i \in \{0,..,T-1\} $ we have that
	\begin{align*}
	\phi_{i}(x,u,y)=\max\{c_{i}(x,u), y\} \ge \max\{c_{i}(x,u), z\}=\phi_{i}(x,u,z),
	\end{align*}
	where the above inequality follows by considering separately the cases $c_{i}(x,u) \ge y$ and $c_{i}(x,u) < y$.
	
	 Assuming $\{c_t\}_{t=0}^T$ are bounded functions the boundedness property, given in Eqn.~\eqref{property: boundednes}, is clearly satisfied by the representation maps given in Eqn.~\eqref{rep: pointwise max} by induction on $i \in\{0,...,T-1\}$. The semi-continuity property (Eqn.~\eqref{property: semi cts}) follow since the point-wise max function, ie $f(x) = \max_{1 \le i \le n}\{x_i\}$, is Lipschitz continuous and hence upper semi-continuous.  $\blacksquare$
\end{pf}
In the next example we consider multiplicative costs. A special case of this cost function, of the form $J(\mbf u, \mbf x) = \mathbb{E}_{\mbf w} [ \exp(\sum_{t=0}^{T-1} c_t(x(t),u(t),w(t)) + c_T(x(T),w(t)))]:= \int \exp(\sum_{t=0}^{T-1} c_t(x(t),u(t),w(t)) + c_T(x(T),w(t))) p(\mbf w) d \mbf w$, where $p(\mbf w)$ is the probability density function of $\mbf w= (w(0),...,w(T))$, has previously appeared \cite{jacobson1973optimal} \cite{glover1988state}.

\begin{ex}[Multiplicative function]
	Suppose $J:  U^T \times \Pi_{t=0}^T X_t \to \R$ is of the form
	{ \begin{align*}
	& J(\mbf u, \mbf x)  = \mathbb{E}_{\mbf w} [c_T(x(T),w(T)) \Pi_{t=0}^{T-1} c_t(x(t),u(t),w(t))) ]\\
	& := \int_{I_0 \times .. I_T} c_T(x(T),w(T)) \Pi_{t=0}^{T-1} c_t(x(t),u(t),w(t)))\\
	& \quad p_T(x(T),w(T)) \Pi_{t=0}^{T-1} p_t(x(t),u(t),w(t)))  dw(0)...dw(T) ,
	\end{align*} }
		where $\mbf u=(u(0),...,u(T-1))$, $\mbf x =(x(0),...,x(T))$, $\mbf w =(w(0),...,w(T))$, $U \subset \R^m$ and $X_t \subset \R^n$, $I_t \subset \R^k$, $c_t: X_t \times U \times I_t \to \R^+$ for $0 \le t \le T-1$, $c_T:X_T \times I_T \to \R$, and $p_t: X_t \times U \times I_t \to \R^+$, $p_T:X_T \times I_T \to \R$ satisfy $\int_{I_t} p_t(x,u,w) dw = 1$ and $\int_{I_T} p_T(x,w) dw = 1$ for $0 \le t \le T-1$ and any $(x,u) \in X_t \times U$. Then $J$ is a monotonically backward separable function. Moreover, if $\{c_t\}_{t=0}^T$ and $\{p_t\}_{t=0}^T$  are bounded functions, and sets $\{I_t\}_{t=0}^T$ have finite measure, then $J$ is naturally monotonically backward separable. Furthermore, if $\int_{I_i}  p_i(x,u,w) c_i(x,u,w) dw \ne 0$ for all $(x,u,i) \in X_i \times U \times \{0,...,T-1\}$ then the associated representation maps are strictly monotonic (Eqn.~\eqref{strict monotonic}).
	
\end{ex}
\begin{pf}
We can write $J$ in Form \eqref{eqn:monotonic backward seperable} using the representation functions

\vspace{-1cm}
	\begin{align} \label{rep: multiplicative}
		\phi_T(x) &  = \int_{I_T} c_T(x)p_T(x,w) dw,\\ \nonumber
		  \phi_{i}(x,u,z) & = \int_{I_i} z p_i(x,u,w) c_i(x,u,w) dw \text{ for } i \in \{0,..,T-1\}.
	\end{align}
	The monotonicity property (Eqn.~\eqref{eqn: monotonic property}) follows as $c_i(x,u,w) \ge 0$ and $p_i(x,u,w) \ge 0$ for all $(x,u,w) \in \R^n \times \R^m \times \R^k$ and $i \in \{0,...,T-1\}$. Furthermore, if $\int_{I_i}  p_i(x,u,w) c_i(x,u,w) dw \ne 0$ for all $(x,u,i) \in X_i \times U \times \{0,...,T-1\}$, then clearly the representation maps are strictly monotonic (Eqn.~\eqref{strict monotonic}).
	
	Assuming $\{c_t\}_{t=0}^T$ and $\{p_t\}_{t=0}^T$ are bounded functions, and sets $\{I_t\}_{t=0}^T$ have finite measure the representation maps in Eqn.~\eqref{rep: multiplicative} clearly satisfy the boundedness property (Eqn.~\eqref{property: boundednes}) by induction on $i \in \{0,...,T-1\}$. For fixed $i \in \{0,..,T-1\}$ and $(x,u) \in X_i \times U$ it follows $\phi_{i}(x,u,z)=cz$, where $c \in \R^+$ is some constant that depends on $(x,u,i)$, is clearly upper semi continuous (as in Eqn.~\eqref{property: semi cts}).  $\blacksquare$
\end{pf}
\vspace{-0.5cm}
In the next example we consider a function that can be interpreted as the expectation of cumulative stochastically stopped additive costs, where at each time stage, $t \in \{0,...,T-1\}$, a cost $c_t(x(t),u(t))$ is added and there is an independent probability, $p_t(x(t),u(t)) \in [0,1]$, of stopping and incurring no further future costs. For a state and input trajectory, $(\mbf u, \mbf x) \in U^T \times \Pi_{t=0}^T X_t$, let us denote the stopping time by $T(\mbf u , \mbf x)$; it then follows the distribution of this random variable is given as
\begin{align} \nonumber
& \mathbb P (T(\mbf u , \mbf x)  = T)  =p_T(x(T)) \Pi_{i=1}^{T-1}(1 - p_i(x(i),u(i))), \\ \label{dist: stopping time}
& \qquad \text{ and for all } t \in \N,  \\ \nonumber
& \mathbb P (T(\mbf u , \mbf x) = t)  = p_t(x(t),u(t)) \Pi_{i=1}^{t-1}(1 - p_i(x(i),u(i))),
\end{align}
where we slightly abuse notation to write $\Pi_{i=1}^{-1}(1 - p_i(x(i),u(i)))=1$ so $\mathbb P (T(\mbf u , \mbf x)  = 0) =p_0(x(0),u(0))$.

The stopped additive function is then given as

	\vspace{-1cm}
\begin{align} \label{11}
&J(\mbf u , \mbf x)= \mathbb E_{T(\mbf u, \mbf x)} \bigg[ \sum_{t=0}^{\min\{T(\mbf u, \mbf x),T-1\}} c_t(x(t),u(t)) \\ \nonumber
& \quad + \mathds{1}_{\{(\mbf u, \mbf x) \in U^T \times \Pi_{t=0}^T X_t : T(\mbf u, \mbf x)=T\}}(\mbf u, \mbf x) c_T(x(T)) \bigg].
\end{align}
To show the function in Eqn.~\eqref{11} is monotonically backward separable we will assume the probability of the stopping time occurring inside the finite time horizon $\{0,...,T\}$ is one; this gives us the following ``law of total probability`` equation $ \sum_{t=0}^T \mathbb P (T(\mbf u, \mbf x)=t)=1 \text{ for all } (\mbf u, \mbf x) \in U^T \times \Pi_{t=0}^T X_t$, which can be rewritten in terms of its probability density functions as,

	\vspace{-1cm}
\begin{align} \nonumber
& \sum_{t=0}^{T-1}p_{t}(x(t),u(t)) \Pi_{i=1}^{t-1}(1 - p_i(x(i),u(i)))\\ \label{law of total probaility}
& \qquad + p_T(x(T))\Pi_{i=1}^{T-1}(1 - p_i(x(i),u(i))) \equiv 1.
\end{align} 
Note, if $p_T(x(T)) \equiv 1$ then it can be trivially shown that Eqn.~\eqref{law of total probaility} holds for any functions $p_i: X_i \times U \to [0,1]$. 

Assuming Eqn.~\eqref{law of total probaility} holds and using the law of total expectation, conditioning on the probability of each stopping time, it follows

\vspace{-0.75cm}

{\small \begin{align} \label{j}
&J(\mbf u , \mbf x)= \mathbb E_{T(\mbf u, \mbf x)} \bigg[ \sum_{t=0}^{\min\{T(\mbf u, \mbf x),T-1\}} c_t(x(t),u(t)) \\ \nonumber
& \qquad + \mathds{1}_{\{(\mbf u, \mbf x) \in U^T \times \Pi_{t=0}^T X_t : T(\mbf u, \mbf x)=T\}}(\mbf u, \mbf x) c_T(x(T)) \bigg]\\ \nonumber
& = \sum_{t=0}^{T-1}  \bigg( \sum_{s=0}^{t} c_s(x(s),u(s)) \bigg)\mathbb P (T(\mbf u, \mbf x)= t) \\ \nonumber
& \qquad + \bigg( \sum_{s=0}^{T} c_s(x(s),u(s)) + c_T(x(T)) \bigg) \mathbb P (T(\mbf u, \mbf x)= T)\\ \nonumber
& =  \sum_{t=0}^{T-1} \bigg( \sum_{s=0}^t c_s(x(s),u(s)) \bigg)  p_{t}(x(t),u(t)) \Pi_{i=0}^{t-1} (1- p_i(x(i),u(i)))\\ \nonumber
& \hspace{1.5cm} + \bigg( \sum_{t=0}^{T-1} c_t(x(t),u(t))  + c_T(x(T)) \bigg)\\ \nonumber
&  \hspace{1.5cm} \qquad \times  p_T(x(T))  \Pi_{i=0}^{T-1} (1- p_i(x(i),u(i))). 
\end{align} }
We next state and prove that the $J$ given in Eqn.~\eqref{j} is monotonically backward separable.
\begin{ex}[Stochastically stopped additive cost] \label{ex: stopped addative}
Suppose $J:  U^T \times \Pi_{t=0}^T X_t \to \R$ is of the form
\vspace{-0.75cm}

{\small\begin{align} \label{fun: stopped addative cost}
&J(\mbf u, \mbf x) = \\ \nonumber &\sum_{t=1}^{T-1} \bigg( \sum_{s=0}^t c_s(x(s),u(s)) \bigg)  p_{t}(x(t),u(t)) \Pi_{i=0}^{t-1} (1- p_i(x(i),u(i)))\\ \nonumber
& \hspace{1.5cm} + \bigg( \sum_{t=0}^{T-1} c_t(x(t),u(t))  + c_T(x(T)) \bigg) \\
& \hspace{1.5cm} \qquad  \times p_T(x(T))  \nonumber \Pi_{i=0}^{T-1} (1- p_i(x(i),u(i)))   , 
\end{align} }
{ \noindent where $p_k: X_k \times U \to [0,1]$ and $p_T: X_T \to [0,1]$ satisfy Eqn.~\eqref{law of total probaility}, $\mbf u=(u(0),...,u(T-1))$, $\mbf x =(x(0),...,x(T))$, $U \subset \R^m$ and $X_t \subset \R^n$, $c_k: X_k \times U \to \R$ and $c_T: X_T \to \R$. Then $J$ is a monotonically backward separable function. Moreover, if $\{c_t\}_{t=0}^T$ are bounded functions, then $J$ is naturally monotonically backward separable. Furthermore, if $p_i(x,u)\ne 1$ for all $(x,u,i) \in X_i \times U \times \{0,...,T-1\}$ then the associated representation maps are strictly monotonic (Eqn.~\eqref{strict monotonic}). } \end{ex}
\vspace{-0.5cm}
\begin{pf}
	Before writing $J$ in the backward separable form (Eqn.~\eqref{eqn:monotonic backward seperable}) we first simplify $J$ by switching the order of the double summation in Eqn.~\eqref{fun: stopped addative cost}. Let $T(\mbf u, \mbf x)$ be a random variable with distribution given in Eqn.~\eqref{dist: stopping time}. As it is assumed $\{p_t\}_{0 \le t \le T}$ satisfy Eqn.~\eqref{law of total probaility} and each time-stage has independent probability of stopping it follows $\sum_{t=s}^{T}\mathbb{P}(T(\mbf u, \mbf x)=t)= \mathbb{P}(T(\mbf u, \mbf x) \ge s )= \mathbb{P}(\cap_{i=0}^{s-1} T(\mbf u, \mbf x) \ne s ) =\Pi_{i=0}^{s-1} \mathbb{P}( T(\mbf u, \mbf x) \ne s ) $. Now,
	
	\vspace{-0.75cm}
	
	{ \small \begin{align*}
	&J(\mbf u, \mbf x)= \sum_{t=0}^{T-1}  \bigg( \sum_{s=0}^{T} c_s(x(s),u(s)) \bigg)\mathbb P (T(\mbf u, \mbf x)= t) \\
	& \qquad \quad + \bigg( \sum_{s=0}^{t} c_s(x(s),u(s)) + c_T(x(T)) \bigg) \mathbb P (T(\mbf u, \mbf x)= T)\\
	&= \sum_{s=0}^{T-1} c_s(x(s),u(s)) \mathbb P(T(\mbf u, \mbf x) \ge s) + c_T(x(T)) \mathbb{P} (T(\mbf u, \mbf x) \ge T)\\
	& = \sum_{s=0}^{T-1} c_s(x(s),u(s)) \Pi_{i=0}^{s-1} \mathbb P( T(\mbf u, \mbf x) \ne i  ) \\
	& \hspace{3.5cm} + c_T(x(T)) \Pi_{i=0}^{T-1} \mathbb P( T(\mbf u, \mbf x) \ne i  ) \\
	& = \sum_{s=0}^{T-1} c_s(x(s),u(s)) \Pi_{i=0}^{s-1} (1-p_i(x(i),u(i)))\\
	& \hspace{2cm}  + c_T(x(T)) p_T(x(T))\Pi_{i=0}^{T-1} (1-p_i(x(i),u(i))).
	\end{align*} }

\vspace{-0.75cm}
	
		It then follows $J$ satisfies Eqn.~\eqref{eqn:monotonic backward seperable} using the representation maps
	\begin{align} \nonumber
	\phi_{i}(x,u,z) & = c_i(x,u) + z (1-p_i(x,u) ) \text{ for } i \in \{0,..,T-1\},\\ \label{phi: min time}
	\phi_T(x) & = c_T(x)p_T(x).
	\end{align}
		The monotonicity property in Eqn.~\eqref{eqn: monotonic property} follows as $(1- p_i(x,u)) \ge 0$ for all $(x,u) \in X_i \times U$ and $i \in \{0,...,T-1\}$. Strict monotonicity (Eqn.~\eqref{strict monotonic}) trivially follows when $p_i(x,u)\ne 1$ for all $(x,u,i) \in X_i \times U \times \{0,...,T-1\}$.
		
		Assuming $\{c_t\}_{t=0}^T$ are bounded functions the representation maps, given in Eqn.~\eqref{phi: min time}, clearly satisfy the boundedness property (Eqn.~\eqref{property: boundednes}) by induction on $i \in\{0,...,T-1\}$. For fixed $i \in \{0,..,T-1\}$ and $(x,u) \in X_i \times U$ it follows $\phi_{i}(x,u,z)=c_0 + c_1z$, where $c_0,c_1 \in \R$ are constants that depends on $(x,u,i)$, clearly satisfies the upper semi continuity property (Eqn.~\eqref{property: semi cts}). $\blacksquare$
	\end{pf}
In the next example we introduce a function representing the number of time-steps a trajectory spends outside some target set. Later, in Section \ref{sec: path planning}, we will use this function as the cost function for path planning problems.

\begin{ex}[Minimum time set entry function] \label{ex: min time}
Suppose $J:  U^T \times \Pi_{t=0}^T X_t \to \R$ is of the form
	\begin{align} \label{fun: min time to get into set}
	J(\mbf u, \mbf x)= \min \bigg\{ \inf \bigg\{t \in [0,T] : x(t) \in S \bigg\}, T \bigg\},
	\end{align}
		where $\mbf u=(u(0),...,u(T-1))$, $u(t) \in \R^m$, $\mbf x =(x(0),...,x(T))$, $x(t) \in \R^n$, $U \subset \R^m$ and $X_t \subset \R^n$, and $S \subset \R^n$. If the set $\{t \in [0,T] : x(t) \in S\}$ is empty, we define the infimum to be infinity. Then $J$ is a naturally monotonically backward separable function.
\end{ex}
%
%
%
%

\begin{pf} 	
	The function given in Eqn.~\eqref{fun: min time to get into set} is actually a special case of the function given in Eqn.~\eqref{fun: stopped addative cost} with  
	\begin{align*}
	p_T(x) & \equiv 1, \quad p_t(x,u)  = \mathds{1}_{S}(x) \text{ for } t\in \{0,...,T-1\}, \\
	c_T(x) & \equiv T, \quad c_t(x,u)  \equiv t.
	\end{align*}
	Note, the functions $\{p_k\}_{0 \le k \le T }$ trivially satisfy Eqn.~\eqref{law of total probaility} as $p_T(x) \equiv 1$. Moreover clearly $\{c_t\}_{t=0}^T$ are bounded functions. Therefore $J$ is naturally monotonically backward separable by Example~\ref{ex: stopped addative}.  $\blacksquare$
\end{pf}

\section{The Principle Of Optimality: A Necessary Condition For Monotonic Backward Separability} \label{sec: POP}
Given a function, $J: \R^{m \times T} \times \R^{n \times (T+1 )} \to \R$, there is no obvious way to determine whether $J$ is monotonically backward separable. Instead, in this section we will recall a necessary condition proposed by Richard Bellman~\cite{bellman1966dynamic}, called the Principle of Optimality (Defn. \ref{defn: principle of optimality deterministic}), that we show all MSOP's with monotonically backward separable cost functions satisfy (Prop.~\ref{prop: P(t,x) satisfies the principle of optimality}). Before recalling the definition of the Principle of Optimality let us consider a family of MSOP's, associated with the tuples $\{J_{t_0}, f , \{X_t\}_{t_0 \le t \le T} , U ,T\}_{t_0=0}^{T}$, each initialized at $(x_0,t_0) \in \R^n \times \{0,....,T\}$, and of the form:
\begin{align}  \label{opt: POP}
&(\mathbf{u}^*,\mathbf{x}^*) {\in} \arg \min_{\mathbf u, \mathbf x} J_{t_0}(\mbf u, \mbf x)  \text{   subject to: }  \\ \nonumber
& x(t+1)=f(x(t),u(t),t) \text{ for  } t={t_0},..,T-1 \\ \nonumber
& x(t_0)=x_0 , \text{ } x(t) \in X_t \subset \mathbb{R}^n \text{ for  } t={t_0},..,T \\ \nonumber
&u(t) \in U \subset \mathbb{R}^m \text{ for  } t={t_0},..,T-1\\ \nonumber
&\mbf u=(u(t_0),...,u(T-1)) \text{ and } \mbf x =(x(t_0),...,x(T))
\end{align}

\begin{defn} \label{defn: principle of optimality deterministic}
	We say the family of MSOP's of Form~\eqref{opt: POP} satisfies \textbf{the Principle of Optimality} at $x_0 \in X_0$ if the following holds. For any $t$ with $0 \le t <T$, if $\mathbf u=(u(0),...,u(T-1))$ and $ \mathbf x=(x(0),...,x(T))$ solve the MSOP initialized at $(x_0,0)$ then $\mathbf v=(u(t),...,u(T-1))$ and $ \mathbf h=(x(t),...,x(T))$ solve the MSOP initialized at $(x(t),t)$.
\end{defn}

\begin{prop} \label{prop: P(t,x) satisfies the principle of optimality} Consider a family of MSOP's of Form~\eqref{opt: POP} associated with $\{ J_{t} , f , \{X_t\}_{t \le s \le T} , U ,T\}_{t=0}^T$. Suppose the MSOP's initialized at $(x_0,0)$ has a unique solution and $J_{t}: U^{T -t} \times \Pi_{s=t}^T X_s \to \R$ is monotonically backward separable (Defn.~\ref{defn: montonic backward seperbale function}). Then the family of MSOP's of Form~\eqref{opt: POP} associated with $\{ J_{t} , f , \{X_t\}_{t \le s \le T} , U ,T\}_{t=0}^T$ satisfies the Principle of Optimality at $x_0 \in X_0$.
\end{prop}
\begin{pf} Since $J_t$ is monotonically backward separable there exists representation maps $\{\phi_t\}_{0 \le t \le T}$ such that 
	
	\vspace{-0.75cm}
	{\small	\[ J_{t}( \mbf u , \mbf x)=\phi_{t}(x({t}),u({t}), \phi_{t+1}(x({t}+1),u({t}+1), \dots \phi_T(x(T)) \dots )). 
		\] } \normalsize
	Now, suppose $\mathbf u^*=(u(0),...,u(T-1))$ and $ \mathbf x^*=(x(0),...,x(T))$ solve the MSOP initialized at $(x_0,0)$ of Form~\eqref{opt: POP} associated with $\{ J_{t} , f , \{X_t\}_{t \le s \le T} , U ,T\}_{t=0}^T$. Suppose for contradiction that there exists some $t \ge 0$ such that $0 \le t  <T$ and $\mathbf v=(u(t),...,u(T-1))$ and $ \mathbf h=(x(t),...,x(T))$ do not solve MSOP initialized at $(x(t),t)$. We will show that this implies that the MSOP initialized at $(x_0,0)$ does not have a unique solution, thus providing a contradiction and verifying the conditions of the Principle of Optimality.  If $(\mathbf v, \mathbf h)$ do not solve MSOP initialized at $(x(t),t)$, then there exist feasible $\mathbf{w}=(w(t),...,w(T-1))$ and $\mathbf z=(z(t),...,z(T))$ such that
	$J_{t}(\mathbf w, \mathbf z) < J_{t}(\mathbf v, \mathbf h)  $. i.e.
	\begin{align} \label{ineq: Js}
	&J_{t}(\mathbf w,\mathbf z )  \\ \nonumber
	& = \phi_{t}(z({t}),w({t}), \phi_{t+1}(z({t}+1),w({t}+1), \dots \phi_T(z(T)) \dots ))\\ \nonumber
	& < \phi_{t}(x({t}),u({t}), \phi_{t+1}(x({t}+1),u({t}+1), \dots \phi_T(x(T)) \dots )) \\ \nonumber
	& =J_{t}(\mathbf v,\mathbf h ).
	\end{align}
	Now, consider the proposed feasible sequences $\hat{\mathbf u}=(u(0),...,u(t-1), w(t),...,w(T-1))$ and $\hat{\mathbf x}=(x(0),...,x(t-1), z(t),...,z(T-1))$. It follows using the monotonicity property (Eqn.~\eqref{eqn: monotonic property}) of monotonically backward separable functions and Inequality~\eqref{ineq: Js}, that
	\begin{align*}
	& J_{0}(\hat{\mathbf u},\hat{\mathbf x} )=\phi_{0}(x({0}),u({0}),\phi_{1}(x(1),u(1),\\
	& \hspace{2cm} \dots \phi_{t}(z({t}),w({t}) \dots \phi_T(z(T)) \dots )) \dots )\\
	&= \phi_{0}(x({0}),u({0}), \dots \phi_{t-1}(x({t-1}),u({t}-1), J_{t}(\mathbf w,\mathbf z ) ) \dots )\\
	& \le \phi_{0}(x({0}),u({0}), \dots \phi_{t-1}(x({t-1}),u({t}-1), J_{t}(\mathbf v,\mathbf h ) ) \dots )\\
	& = J_{0}(\mathbf u^*,\mathbf x^* ),
	\end{align*}
	which shows $(\hat{\mathbf u},\hat{\mathbf x} )$ is optimal  contradicting that $(\mbf u^*, \mbf x^*)$ is the unique solution of the MSOP at $(x_0,0)$.  $\blacksquare$
\end{pf} 
Prop.~\ref{prop: P(t,x) satisfies the principle of optimality} shows the Principle of Optimality (Defn.~\ref{defn: principle of optimality deterministic}) is a necessary condition that all families of MSOP's with unique solutions and monotonically backward separable cost functions must satisfy. We now conjecture a necessary and sufficient condition. The following notation is used in this conjecture. Given $J_t$, $\{X_t\}_{0 \le t \le T}$ and $U$ let us denote the set $\mcl F$, where $(f,x_0) \in \mcl F$ if $x_0 \in X_0$ and the MSOP associated with $\{J_{0}, f , \{X_t\}_{0 \le t \le T} , U ,T\}$ initialized at $(x_0,0)$ has a unique solution.
\begin{conj} \label{conj: MBS iff PO}
	Consider $\{X_t\}_{0 \le t \le T} \subset \R^{n \times T}$, $U \subset \R^m$ and $J_{t}: U^{T -t} \times \Pi_{s=t}^T X_s \to \R$. Then, for any $(f,x_0) \in \mcl F$ the family of MSOP's associated with $\{J_{t}, f , \{X_t\}_{t \le s \le T} , U ,T\}_{t=0}^{T}$ satisfy the Principle of Optimality at $x_0 \in X_0$ if and only if $J_t$ is monotonically backward separable.
	\end{conj}
Regardless of whether Conjecture~\ref{conj: MBS iff PO} is true, Prop.~\ref{prop: P(t,x) satisfies the principle of optimality} is useful. Prop.~\ref{prop: P(t,x) satisfies the principle of optimality} provides a way of proving a function $J_{t}:  U^{T -t} \times \Pi_{s=t}^T X_s \to \R$ is not monotonically backward separable. Rather than showing $J_{t}$ does not satisfy Defn.~\ref{defn: montonic backward seperbale function} for every family of representation maps $\{\phi_s \}_{s=t}^T$, for which there are an uncountably many, we find any $f$ for which the family of MSOP's associated with $\{J_{t}, f , \{X_s\}_{t \le s \le T} , U ,T\}_{t=0}^{T}$ has a unique solution for some initialization $(x_0,0)$ and does not satisfy the Principle of Optimality. Then Prop.~\ref{prop: P(t,x) satisfies the principle of optimality} shows $J_{t}$ is not monotonically backward separable. We demonstrate this proof strategy in the following lemma.
\begin{lem} \label{lem: functions that do not satisfy PO}
	The function $J_{t}: \R^{m \times (T -t)} \times \R^{n \times (T+1 -t )} \to \R$, defined as
	
	\vspace{-1cm}
	\begin{align} \label{fun: non MBS}
	J_{t}(\mbf u, \mbf x): = \sum_{s= t}^{T-1} c_s(u(s)) + \max_{t \le s \le T} d(x(s)),
	\end{align}
	is not monotonically backward separable (Defn.~\ref{defn: montonic backward seperbale function}) for all functions $c_k: \R^m \to \R$ and $d_k: \R^n \to \R$.
	\end{lem}
\begin{pf}
	Let $T=3$, $n=1$ and $m=1$. Consider the cost functions $c_{0}(u)=-u,$ $c_{1}(u)=u,$ $c_{2}(u)=-u/2$, and $d(x)=x$. Consider the dynamics $f(x,u,t)=x+u$ and constraints $X_t=[0,h]$ and $U=\{-h,0,h\}$, where $h>0$. Let us consider the MSOP of Form~\eqref{opt: POP} associated with $\{J_0, f , \{X_t\}_{0 \le t \le 3} , U ,3\}$ and initialized at $(x_0,t_0)=(0,0)$. It can be shown there are $3^3=27$ input sequences in $\{-h,0,h\}^3$, only 8 of which are feasible to the MSOP initialized at $(x_0,t_0)=(0,0)$. By calculating the cost of each feasible input we can deduce the unique optimal input sequence is $\mbf u= (h,-h,h)$, yielding a unique optimal trajectory of $\mbf x=(0,h,0,h)$. Following the input sequence to $t=2$ we examine the MSOP of Form~\eqref{opt: POP} initialized at $(x_0,t_0)=(0,2)$. For the MSOP initialized at $(x_0,t_0)=(0,2)$ there are only two feasible inputs: $u(2)=0$ or $u(2)=h$. Of these, the first is optimal (cost of 0 vs $h/{2}$). Thus although $\mbf u= (h,-h,h)$ and $\mbf x=(0,h,0,h)$ solve the MSOP initialized at $(x_0,t_0)=(0,0)$, $\mbf v = (h)$ and $\mbf h= (0,h)$ do not solve the MSOP initialized at $(x_0,t_0)=(0,2)$. We conclude the family of MSOP's associated with $\{J_{t}, f , \{X_s\}_{t \le s \le 3} , U ,3\}_{t=0}^{3}$ does not satisfy the Principle of Optimality at $x_0=0$, although the MSOP initialized at $(x_0,t_0)=(0,0)$ does have a unique solution. Therefore by Prop.~\ref{prop: P(t,x) satisfies the principle of optimality} the function $J_{t}$ is not monotonically backward separable.  $\blacksquare$	\end{pf}

\begin{rem}
	The function given in Eqn.~\eqref{fun: non MBS} can clearly be expressed as the addition of two monotonically backward separable functions, $J_1(\mbf u, \mbf x)=\sum_{s= t}^{T-1} c_s(u(s))$ (Lemma~\ref{lem: add sep is back sep}) and $J_2(\mbf u, \mbf x)= \max_{t \le s \le T} d(x(s))$ (Example~\ref{ex: sup forward sep}). Therefore, Lemma~\ref{lem: functions that do not satisfy PO} shows that the property of monotonically backward separability is not preserved under addition.
	\end{rem}
\section{Comparison With State Augmentation Methods} \label{sec: compare}
We proposed an alternative method for solving MSOP's with non-additively separable costs in \cite{jones2018extensions}; where cost functions are forward separable:
  	\begin{align} \nonumber
&J(\mbf u, \mbf x)=\psi_T( x(T),\psi_{T-1}(x(T-1),u(T-1),\psi_{T-2}(....,\\ \label{fun: forward sep}
&\qquad  \psi_{1}(x(1),u(1), \psi_{0}(x(0),u(0)))),....,))), 
\end{align}
where $\psi_0 : X_0 \times U \to \R^k$, $\psi_t : X_{t} \times U \times Image\{\psi_{t-1}\} \to \R^k$ for $ t \in \{1,..,T-1\}$, and $\psi_T: X_{T} \times Image\{\psi_{T-1}\}   \to \R$.

It was shown that for $\{J, f , \{X_t\}_{0 \le t \le T} , U ,T\}$, where $J$ is of the Form~\eqref{fun: forward sep}, an equivalent MSOP with additively separable cost function, $\{\tilde{J}, \tilde{f} , \{\tilde{X}_t\}_{0 \le t \le T} , U ,T\}$, can be constructed, where $\tilde{J}(\mbf u, \mbf x)= \psi_T(x(T))$,  $\tilde{f}([x_1,x_2]^T,u,t)= [f(x_1,u,t), \psi_{t}(x_1,u,x_2) ]^T$, and $\tilde{X}_t= X_t \times Image\{\psi_{t}\}$. The augmented MSOP, $\{\tilde{J}, \tilde{f} , \{\tilde{X}_t\}_{0 \le t \le T} , U ,T\}$, can then be solved using the classical Bellman Equation \eqref{eqn: Bellman Eqn}. This state augmentation method is particularly useful when solving MSOP's with cost functions that are not monotonically backward separable, for instance the function in Eqn.~\eqref{fun: non MBS}. However, the augmented MSOP has a larger state space dimension. Therefore, in the case when the cost function is both forward separable, of Form \eqref{fun: forward sep}, and monotonically backward separable, of Form \eqref{eqn:monotonic backward seperable}, it is computationally more efficient to solve the GBE \eqref{eqn: generalize Bellman} rather than augmenting and solving Bellman's Equation \eqref{eqn: Bellman Eqn}. We now demonstrate this in the following numerical example.

\begin{figure} 	
	\centering
	\includegraphics[scale=0.6]{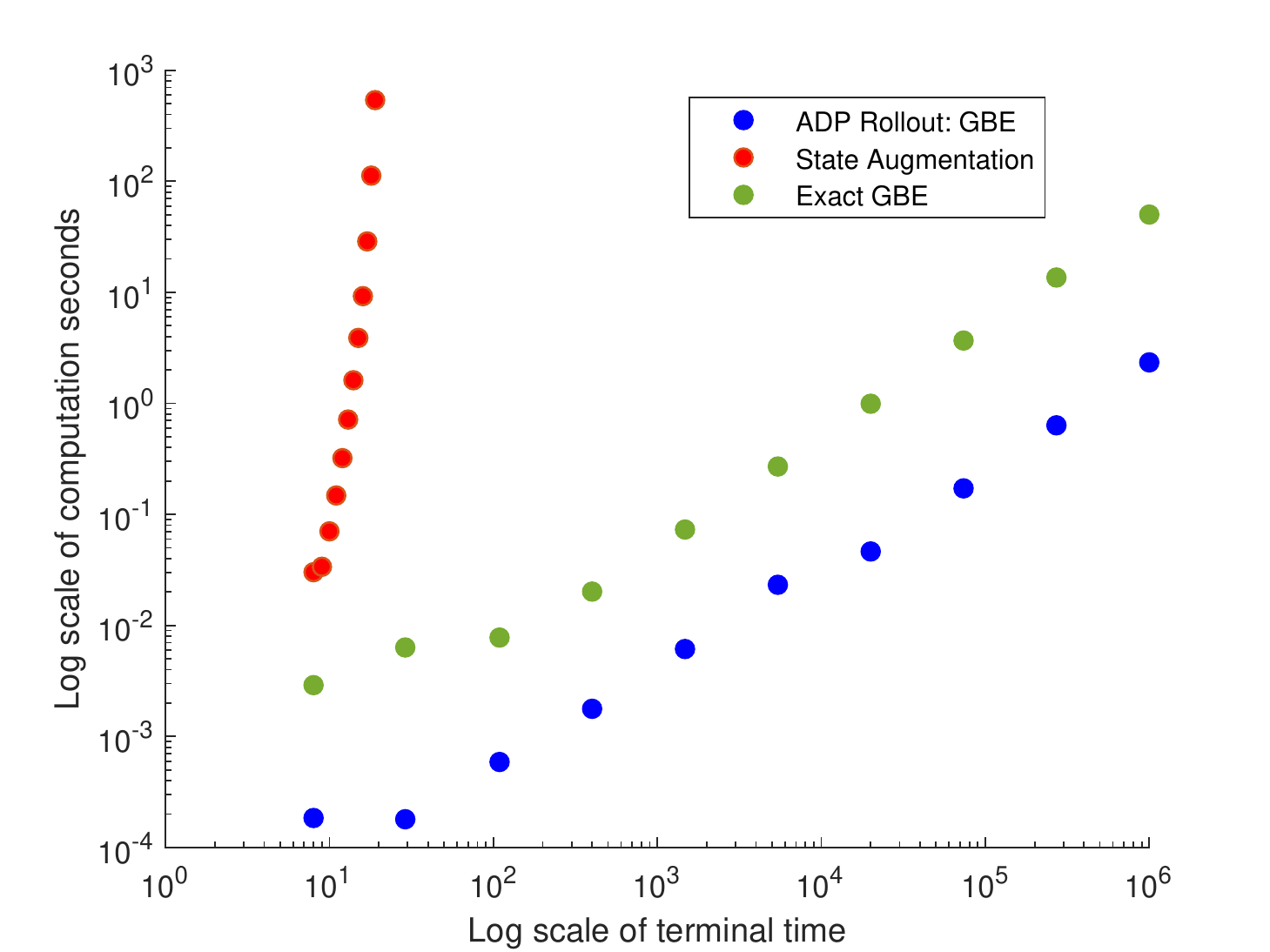}
	\vspace{-10pt}
	\caption{Log log graph showing computation time for solving MSOP \eqref{MSOP:sqrt} using state augmentation (red points), via exactly solving GBE (green points), and via approximately solving the GBE using the rollout (blue points) algorithm versus the terminal time of the problem.}
	\label{fig: state_aug_vs_GBE}
	\vspace{-5pt}
\end{figure}

Consider the MSOP 
{ \begin{align} \nonumber
& \min_{\mathbf u, \mathbf x} \sqrt{x(0) + u(0) + \sqrt{..... \sqrt{x(T-1) + u(T-1) + \sqrt{x(T) }}}} \\ \label{MSOP:sqrt}
&\text{   subject to: }  \\ \nonumber
& x(t+1)=\begin{cases}
2 \text{ if } u = 0.5\\
1 \text{ if } u= 1
\end{cases} \text{ for  } t={0},..,T, \\ \nonumber
& x(0)=2 , \text{ } x(t) \in \{1,2\}  \text{ for  } t={0},..,T, \\ \nonumber
&u(t) \in \{0.5,1\} \text{ for  } t={0},..,T-1.
\end{align} }
The cost function in the above MSOP is naturally monotonically backward separable and can be written in the Form \eqref{eqn:monotonic backward seperable} with representation maps
\begin{align} \label{fun: rep sqrt}
\phi_T(x)= \sqrt{x}, \text{ } \phi_t(x,u,z)= \sqrt{x + u + z} \text{ for } t \in \{0,..,T-1\}.
\end{align}
Moreover the cost function is also forward separable and can be written in the Form \eqref{fun: forward sep} with representation maps
\begin{align} \label{fun: state aug sqrt}
\psi_0(x,u)=[x, u]^T, \quad \psi_t(x,u,z)=[z,x,u]^T, \\ \nonumber
\psi_T(x,z)= \sqrt{z_1 + z_2 + \sqrt{.... \sqrt{z_{2T-1} + z_{2T} + \sqrt{x}}}}.
\end{align}
We solved the MSOP in Eqn.~\eqref{MSOP:sqrt} using both the GBE and the state augmentation method, plotting the computation time results in Figure~\ref{fig: state_aug_vs_GBE}. The green points represent the computation time required to construct the value function by solving the GBE~\eqref{eqn: generalize Bellman} with representation maps given in Eqn.~\eqref{fun: rep sqrt}, and then to synthesize the optimal input sequence using Eqn~\eqref{fun: optimal policy}. The red points represent the computation time required to construct the value function by solving Bellman's Equation~\eqref{eqn: Bellman Eqn} for the state augmented MSOP and then to construct the optimal input sequence. The green points increases linearly as a function of the terminal time, $T \in \N$, of order $\mathcal{O}(T)$, whereas the red points increases exponentially with respect to $T$, of order $\mathcal{O}(2^{T})$ (due to the fact that using representation maps, given in Eqn.~\eqref{fun: state aug sqrt}, results in an augmented state space of size $2^T$). Moreover, Figure~\ref{fig: state_aug_vs_GBE} also includes blue dots representing computation times required to solve the GBE approximately, as discussed in the next section.

\subsection{Approximate Dynamic Programming Using The GBE}
Rather than solving the MSOP~\eqref{MSOP:sqrt} exactly using the GBE, as we did in the previous section, we now use an Approximate Dynamic Programming (ADP)/Reinforcement Learning (RL) algorithm to heuristically solve the MSOP and numerically show these algorithms can result in lower computational times when compared to methods that solve the GBE exactly. This demonstrates that MSOP's with monotonically backward separable cost functions can be heuristically solved using the same methods developed in the ADP literature with the aid of the methodology developed in this paper.

 Typically ADP methods use parametric function fitting (neural networks, linear combinations of basis functions, decision tree's, etc) to approximate the value function from data. The approximated value function is then used to synthesize a sub-optimal input sequence. To see how this works, suppose an ADP algorithm constructs some approximate value function, denoted $\tilde{V}(x,t)$, then an approximate optimal input sequence, $\mbf{\tilde{u}}=(\tilde{u}(0),...,\tilde{u}(T))$, can be constructed by solving
		\begin{align} \nonumber 
& \tilde{u}(k) \in   \arg \inf_{ u \in \Gamma_{\tilde{x}(k),k}} \bigg\{ \phi_t(\tilde{x}(k),u,\tilde{V}(f(\tilde{x}(k),u,k),k+1)) \bigg\}\\ \nonumber
& \hspace {4cm} \text{ for } k\in \{0,...,T-1\}.\\
\nonumber & \tilde{x}(0)=x_0, \quad \tilde{x}(k+1)= f(\tilde{x}(k),\tilde{u}(k),k) \label{rollout policy} \\
& \hspace {4cm} \text{ for } k\in \{0,...,T-1\}.
\end{align}
One way to obtain an approximate value function, $\tilde{V}$, is to use the rollout algorithm found in the textbook \cite{bertsekas1995dynamic}. This algorithm supposes a base policy is known $\mu_{base}: \R^n \times \N \to U$ and approximates the value function as follows
\begin{align*}
 \tilde{V}(x,t) = &  \phi_t(  x(t),u(t),  \phi_{t+1}(x(t+1),u(t+1),... \phi_T(x(T))...)),\\
&\text{where } x(t)=x \text{ and for all }  s \in \{t,...,T-1\}, \\
 & x(s+1)=f(x(s),u(s),t), \text{ } u(s) = \mu_{base}(x(s),s).
\end{align*}

Using the base policy $\mu_{base}(x,t)=\begin{cases}
1 \text{ if } t/4 \in \N\\
0.5 \text{ otherwise}
\end{cases} $ we used the rollout algorithm to solve the MSOP~\eqref{MSOP:sqrt} for terminal times $T=8$ to $10^{6}$. Computation times are plotted as the blue points in Figure~\ref{fig: state_aug_vs_GBE} showing better performance than solving the GBE exactly or using state augmentation.

\section{Application: Path Planning And Obstacle Avoidance} \label{sec: path planning}
In this section we design a full state feedback controller (Markov Policy) for a discrete time dynamical system with the objective of reaching a target set in minimum time while avoiding moving obstacles.
\subsection{MSOP's For Path Planning}

We say the MSOP, associated with tuple $\{J, f , \{X_t\}_{0 \le t \le T} , U ,T\}$, defines a Path Planning DP problem if
\begin{itemize}
	\item $J(\mbf u, \mbf x)= \min \bigg\{ \inf \bigg\{t \in [0,T] : x(t) \in S \bigg\}, T \bigg\}$.
	\item $S=\{x \in \R^n: g(x) < 0\}$, where $g: \R^n \to \R$.
	\item  $X_t =\R^n/(\cup_{i=1}^N O_{t,i})$, where $O_{t,i}= \{x \in \R^n: h_{t,i}(x)<0 \}$ and $h_{t,i}: \R^n \to \R$.
	\item There exits a feasible solution, $(\mathbf u, \mathbf x)$, to the MSOP \eqref{opt: DP} associated with the tuple $\{J, f , \{X_t\}_{0 \le t \le T} , U ,T\}$ such that $x(k) \in S$ for some $k \in \{0,...,T\}$.
\end{itemize}
Clearly, solving the MSOP~\eqref{opt: DP} associated with a path planning problem tuple, $\{J, f , \{X_t\}_{0 \le t \le T} , U ,T\}$, is equivalent to finding the input sequence that drives a discrete time system, governed by the vector field $f$, to a target $S$ in minimum time while avoiding the moving obstacles, represented as sets  $O_{t,i} \subset \R^n$. Moreover, as shown in Example~\ref{ex: min time}, the cost function $J$ is a naturally forward separable function (Defn.~\ref{defn: montonic backward seperbale function}).

\subsection{Path Planning for Dubin's Car}
We now solve the path planning problem with dynamics as defined in \cite{maidens2018symmetry}; also known as the Dubin's car dynamics.
\begin{align} \label{dyn: dubin}
f\left( x,u,t \right)= 
\left[x_1 + v \cos(x_3), x_2 + v \sin(x_3), x_3 + \frac{v}{L} \tan(u)  \right]^T,
\end{align}
where $(x_1, x_2) \in \R^2$ is the position of the car, $x_3 \in \R$ denotes the angle the car is pointing, $u \in \R$ is the steering angle input, $v \in \R$ is the fixed speed of the car, and $L$ is a parameter that determines the turning radius of the car.

We solve the path planning problem using a discretization scheme, similar to \cite{jones2018extensions}; such discretization schemes are known to be parallelizable \cite{maidens2016parallel}. The target set, obstacles, state space, and input constraint sets are given by
\begin{align*}
S & = \{(x_1,x_2) \in \R^2: -0.25 <x_1-0.75 <0.25,\\
& \hspace{4cm} -0.25<x_2+ 0.75< 0.25 \}\\
O_{t,i} & = \{(x_1,x_2) \in \R^2:  (x_1-X_i)^2 + (x_2- Y_i)^2 -R_i^2<0  \}\\
& \hspace{1cm} \text{for } i \in \{1,...,15\} \text{ and } t \in \{0,...,T\}\\
& X_t=[-1,1]^2 \times \R \text{ for } t \in \{0,...,T\}, \quad U= [-1,1],
\end{align*}
where $X,Y,R\in \R^{15}$ are randomly generated vectors. The parameters of the system are set to $v=0.1$ and $L=1/6$.

	\begin{figure} 	
	\centering
	\includegraphics[scale=0.6]{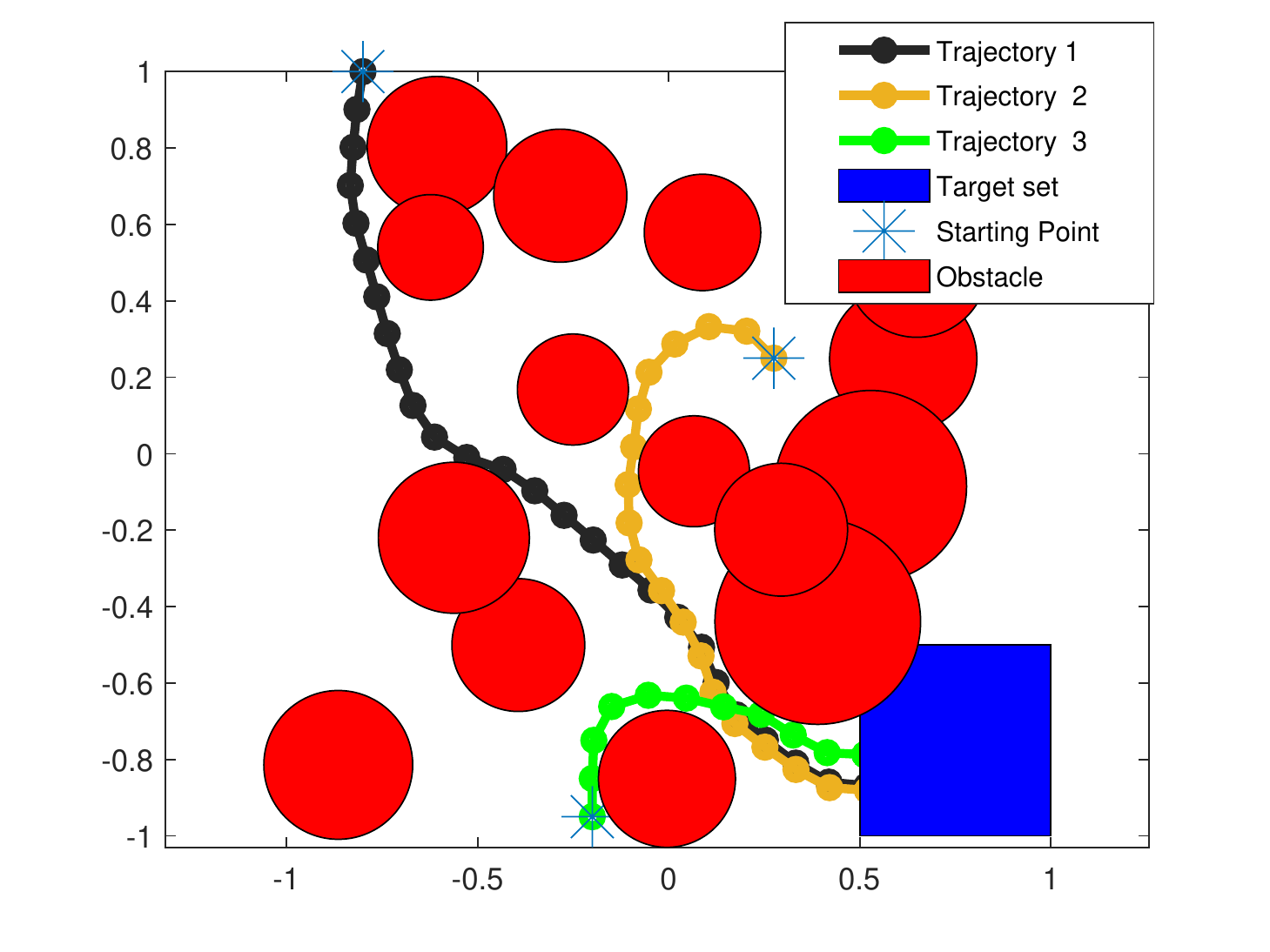}
	\vspace{-10pt}
	\caption{Graph showing approximate optimal trajectories, shown as the gold, black and green curves, with dynamics given in Eqn.~\eqref{dyn: dubin} and the goal of reaching the target set, shown as the blue square, while avoiding obstacles, shown as red circles. }
	\label{fig: dubin car 1}
	\vspace{-5pt}
\end{figure}

	\begin{figure} 	
	\centering
	\includegraphics[scale=0.6]{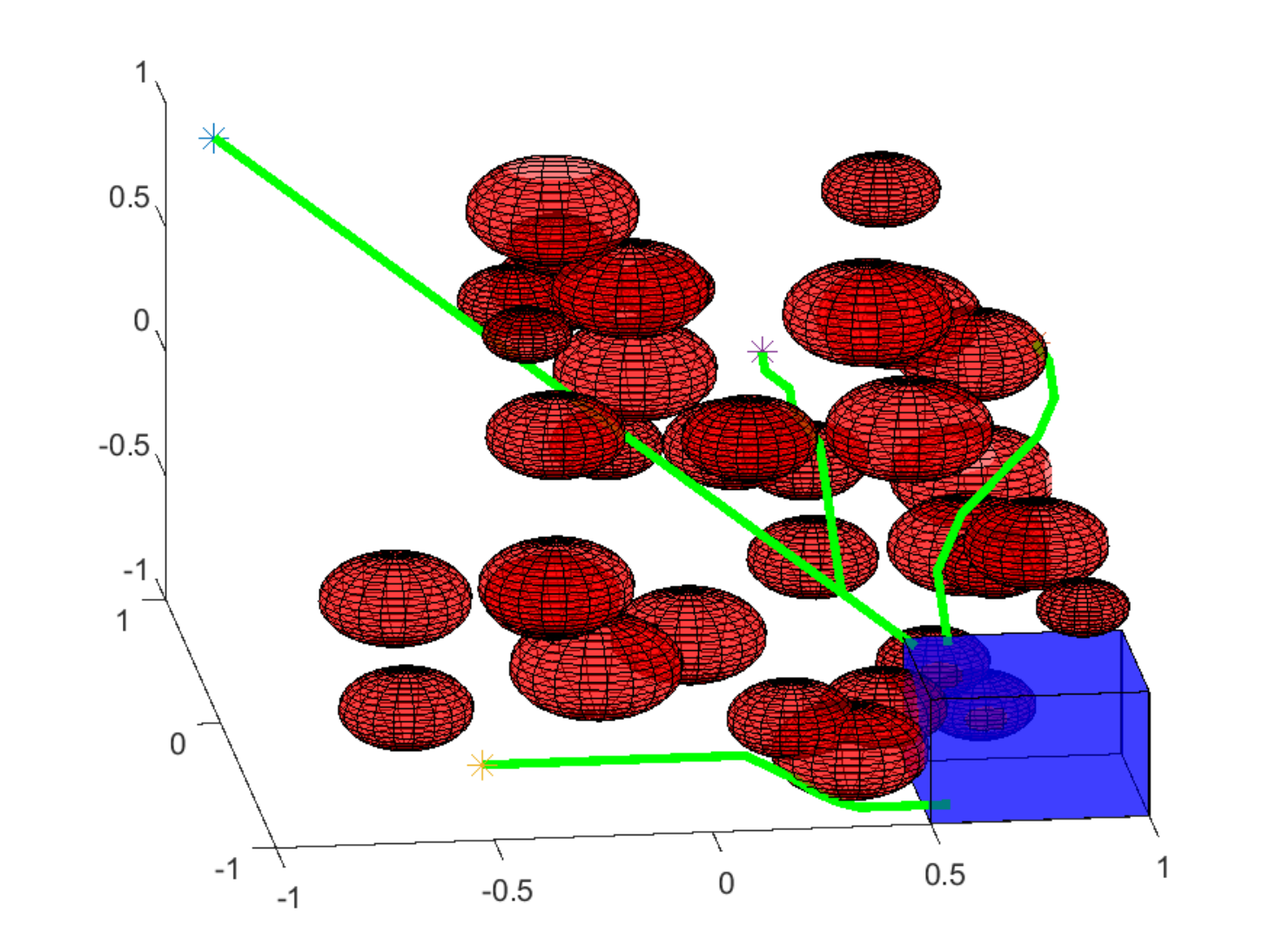}
	\vspace{-20pt}
	\caption{Graph showing approximate optimal trajectories, shown as the green curves, with dynamics given in Eqn.~\eqref{dyn: king dynamics} and the goal of reaching the target set, shown as the blue cube, while avoiding obstacles, shown as red spheres.} \label{fig: 3D king}
	\vspace{-5pt}
\end{figure}

Figure \ref{fig: dubin car 1} shows three approximately optimal state sequences starting from different initial conditions. These state sequences are found by numerically solving the GBE~\eqref{eqn: generalize Bellman}, where $\{\phi_t\}_{t=0}^T$ are as in Example~\ref{ex: min time}. To numerically solve the GBE~\eqref{eqn: generalize Bellman} the state space, $X_t \subset \R^3$, is discretized as a $60 \times 60 \times 60$-grid between $[-1,1]^2 \times [0,2\pi]$ and the input space, $U \subset \R$, is discretized as $100$ grid points within $[-1,1]$. The first state sequence was chosen to have initial condition $[-0.8,1,-0.55 \pi]^T \in \R^3$ (the furthest of the three trajectories from the target) and took $25$ steps to reach its goal. The second state sequence was chosen to have initial condition $[0.275,0.25,0.75 \pi]^T \in \R^3$; in this case as $x_3(0)=0.75 \pi$ Dunbin's car initially is directed towards the top left corner. The input sequence successfully turns the car downwards between two obstacles and into the target set, taking $18$ steps. The third trajectory was chosen to have initial condition $[-0.2,0.95,0.5 \pi]^T \in \R^3$-starting very closely to an obstacle facing upwards. This trajectory had to use the full turning radius of the car to navigate around the obstacle towards the target set and took $10$ steps.

\subsection{Path Planning in 3D}

We now solve a three dimensional path planning problem with dynamics given by
\begin{align} \label{dyn: king dynamics}
f\left( x,u,t \right)= 
[x_1 + u_1, x_2 + u_2,x_3 + u_3  ]^T.
\end{align}
The target set, obstacles, state space and input constraint set were respectively are given by
\begin{align*}
&S  = \{(x_1,x_2,x_3) \in \R^2: -0.25 <x_1-0.75 <0.25,\\
& \quad -0.25<x_2+ 0.75< 0.25, -0.25<x_2+ 0.75< 0.25 \}\\
&O_{t,i}  = \{(x_1,x_2,x_3) \in \R^3:  (x_1-A_i -\alpha_i t )^2 + (x_2- B_i -\beta_i t)^2 \\
& + (x_2- C_i -\gamma_i t)^2-R_i^2<0  \} \hspace{-0.1cm} \text{ for } \hspace{-0.05cm} i \in \{1,...,35\}, t \in \{0,...,T\} \\
& X_t-[-1,1]^3 \text{ for } t \in \{0,...,T\}, \qquad U= [-0.05,0.05]^3,
\end{align*}
where $A,B,C,\alpha,\beta,\gamma,R\in \R^{35}$ are randomly generated vectors. Note, when $\alpha,\beta,\gamma$ are non-zero the center of the spherical obstacles moves with time. For presentation purposes in this paper we consider stationary obstacles, selecting $\alpha=\beta=\gamma=0$, however, a downloadable .gif file showing the numerical solution for moving obstacles can be found at~\cite{data}.

This path planning problem can be numerically solved by computing the solution to the GBE~\eqref{eqn: generalize Bellman} using $\{\phi_t\}_{t=0}^T$ as given in Example~\ref{ex: min time}. To numerically solve the GBE~\eqref{eqn: generalize Bellman} we discretized the state and input space, $X_t \subset \R$ and $U \subset \R^3$, as a $40 \times 40 \times 40$ uniform grid on $[-1,1]^3$ and a $5 \times 5 \times 5$ uniform grid on $ [-0.05,0.05]^3$ respectively. Figure~\ref{fig: 3D king} shows four optimal state sequences, shown as green lines, starting from various initial conditions. All trajectories successfully avoid the obstacles, represented as red spheres, and reach the target set, shown as a blue cube.

\textbf{GPU Implementation} All DP methods involving discretization fall prey to the curse of dimensionality, where the number of points required to sample a space increases exponentially with respect to the dimension of the space. For this reason solving MSOP's in dimensions greater than three can be computationally challenging. Fortunately, our discretization approach to solving the GBE~\eqref{eqn: generalize Bellman}, can be parallelized at each time-step. To improve the scalability of the proposed approach, we have therefore constructed in Matlab a GPU accelerated DP algorithm for solving the 3D path planning problem. This code is available for download at Code Ocean~\cite{CodeOcean}.

\section{Application: Maximal Invariant Sets} \label{sec: invariant sets}
The Finite Time Horizon Maximal Invariant Set (FTHMIS) is the largest set of initial conditions such that there exists an input sequence that produces a feasible state sequence over a finite time period. Computation of the maximal robust invariant sets over infinite time horizons was considered in \cite{xue2018robust}. Before we define the FTHMIS we introduce some notation.

For  $f: \R^n \times \R^m \times \N \to \R^n$ we say the map $\rho_f: \R^n \times \R \times \R^{m \times (T-1)} \to \R^n$ is the solution map associated with $f$ if the following holds for all $(x_0,t, \mbf u) \in \R^n \times \{0,...,T\} \times \R^{m \times (T-1)} $
\begin{align*}
\rho_f(x_0,t,\mathbf u)= x(t),
\end{align*}
where $\mathbf u=(u(0),...,u(T-1))$, $x(k+1)=f(x(k),u(k),k)$ for all $k \in \{0,..,k-1\}$, and $x(0)=x_0$.

\begin{defn} \label{defn: finite time maximal invariant set}
	For $f: \R^n \times \R^m \times \N \to \R^n$, $X_t \subseteq \R^n$, $U \subset \R^m$, $T \in \N$, and $\mathcal A_t \subseteq \R^n$ we define \textit{the Finite Time Horizon Maximal Invariant Set} (FTHMIS), denoted by $\mathcal{R}$, by
	\begin{align*}
	\mathcal{R}:= & \{x_0 \in \R^n: \text{there exists }\mathbf u \in \Gamma_{x_0,[0,T-1]} \text{ such that } \\
	& \hspace{2cm} \rho_f(x_0,t,\mathbf u) \in \mathcal A_t \text{ for all } t \in \{0,...,T\} \},
	\end{align*}
	where the notation $\Gamma_{x_0,[0,T-1]}$ is as in Eqn.~\eqref{notation: feasible input sequence}.
\end{defn}


We next show that the sublevel set of the value function associated with a certain MSOP can completely characterize the FTHMIS.

\begin{thm} \label{thm: FTHMIS characterization}
	Consider the sets $\mathcal A_t=\{x \in \R^n: g_t(x) < 0 \}$, where $g_t:\R^n \to \R$. Suppose $V$ is a value function associated with the MSOP, defined by the tuple $\{J, f , \{X_t\}_{0 \le t \le T} , U ,T\}$, where $J(\mbf u, \mbf x)=\max_{0 \le k \le T}{g_k(x(k))}$. Then
	
	\vspace{-1cm}
	\begin{align} \label{FTHRMIS characterization}
\mathcal{R}= \{x\in \R^n: V(x,0) < 0\},
	\end{align}
	where the set $\mathcal{R} \subset \R^n$ is the FTHMIS as in Defn.~\ref{defn: finite time maximal invariant set}.
\end{thm}
\begin{pf}
	The function $J(\mbf u, \mbf x)=\max_{0 \le k \le T}{g_k(x(k))}$ is monotonically backward separable as shown in Example~\ref{ex: sup forward sep} using representation maps given by
		\begin{align*}
	\phi_{i}(x,u,z) & =\max\{g_{i}(x), z\} \text{ for all } i \in \{0,..,T-1\} \\
	\phi_T(x) & = g_T(x).
	\end{align*}
	Therefore by Defn.~\ref{defn: Value functions} any value function, $V: \R^n \to \R$, associated with $\{J, f , \{X_t\}_{0 \le t \le T} , U ,T\}$ satisfies
	\begin{align}
	V(x,T)= g_T(x) \text{ for all } x \in X_T,
	\end{align}
	and for all $t \in \{0,1,..,T-1\}$ and $x \in X_t$
	\begin{align}
	&V(x,t)= \inf_{\mathbf u \in \Gamma_{x,[0,T-1]}} \max_{t \le k \le T}{g_k(\rho_f(x,k,\mathbf u))}.
	\end{align}
We will first show that $\mathcal{R} \subseteq \{x\in \R^n: V(x,0)<0\}$. Let $x_0 \in \mathcal{R}$ then by Defn.~\ref{defn: finite time maximal invariant set} there exists $\mathbf{u}_0 \in \Gamma_{x_0,[0,T-1]}$ such that

\vspace{-1cm}
\begin{align*}
\rho_f(x_0,t,\mathbf u_0) \in \mathcal A_t \quad \text{for all } t \in \{0,...,T\}.
\end{align*}
As $\mathcal A_t=\{x \in \R^n: g_t(x) < 0 \}$ we deduce from the above equation that
\begin{align} \label{1}
g_t(\rho_f(x_0,t,\mathbf u_0)) < 0 \quad \text{for all } t \in \{0,...,T\}.
\end{align}
Therefore,

\vspace{-1cm}
\begin{align*}
V(x_0,0) & =\inf_{\mathbf u \in \Gamma_{x_0,[0,T-1]}} \max_{0 \le k \le T}{g_k(\rho_f(x_0,k,\mathbf u))}\\
& \le \max_{0 \le k \le T}{g_k(\rho_f(x_0,k,\mathbf u_0))} < 0,
\end{align*}
where the second inequality follows by Eqn.~\eqref{1}. We therefore deduce $x_0 \in \{x \in \R^n: V(x,0) < 0\}$ and hence $\mathcal{R} \subseteq \{x\in \R^n: V(x,0) < 0\}$.

We next show $\{x\in \R^n: V(x,0) < 0\}  \subseteq  \mathcal{R}$. Let $x_0 \in \{x\in \R^n: V(x,0) < 0\}$ then,
\begin{align*}
	\inf_{\mathbf u \in \Gamma_{x_0,[0,T-1]}} \max_{0 \le k \le T}{g_k(\rho_f(x_0,k,\mathbf u))}= V(x_0,0) < 0.
\end{align*}
Therefore as the above inequality is strict, there exists some $ \eps>0$ such that
\begin{align} \label{ineq: less than eps}
\inf_{\mathbf u \in \Gamma_{x_0}} \max_{0 \le k \le T}{g_k(\rho_f(x_0,k,\mathbf u))}= V(x_0,0) < -\eps.
\end{align}
By the definition of the infimum for any $\delta>0$ there exits $\mathbf{w} \in \Gamma_{x_0,[0,T-1]}$ such that
\vspace{-0.8cm}
{\small\begin{align} \label{inf property}
\max_{0 \le k \le T}{g_k(\rho_f(x_0,k,\mathbf w )) } & < \inf_{\mathbf u \in \Gamma_{x_0,[0,T-1]}} \max_{0 \le k \le T}{g_k(\rho_f(x_0,k,\mathbf u))} + \delta.
\end{align}} 

\vspace{-1cm}
Hence by letting $0<\delta<\eps$ we get
\vspace{-0.8cm}
{\small\begin{align} \nonumber
\max_{0 \le k \le T}{g_k(\rho_f(x_0,k,\mathbf w )) } & < \inf_{\mathbf u \in \Gamma_{x_0,[0,T-1]}} \max_{0 \le k \le T}{g_k(\rho_f(x_0,k,\mathbf u))} + \delta\\
 & <-\eps + \delta  <0, \label{label}
\end{align} }
\vspace{-1.1cm}

{where the first inequality follows by Eqn.~\eqref{inf property}, the second inequality follows from Eqn.~\eqref{ineq: less than eps}, and the third inequality follows from selecting $\delta<\eps$. \noindent}

Therefore by Eqn.~\eqref{label} there exists $\mathbf{w} \in \Gamma_{x_0,[0,T-1]}$ such that $\max_{0 \le k \le T}g_k(\rho_f(x_0,k,\mathbf w)) < 0$. We now deduce that for any $t \in \{0,...,T\}$
\begin{align*}
g_t(\rho_f(x_0,t,\mathbf w)) \le \max_{0 \le k \le T}g_k(\rho_f(x_0,k,\mathbf w)) < 0.
\end{align*}
Thus $\rho_f(x_0,t,\mathbf u_0) \in \mathcal A_t$, implying $x_0 \in \mathcal{R}$. Therefore $\{x\in \R^n: V(x,0) < 0\}  \subseteq  \mathcal{R}$. $\blacksquare$
\end{pf}

\subsection{Numerical Example: Maximal Invariant Sets}
\begin{figure} 	
	\centering
	\includegraphics[scale=0.6]{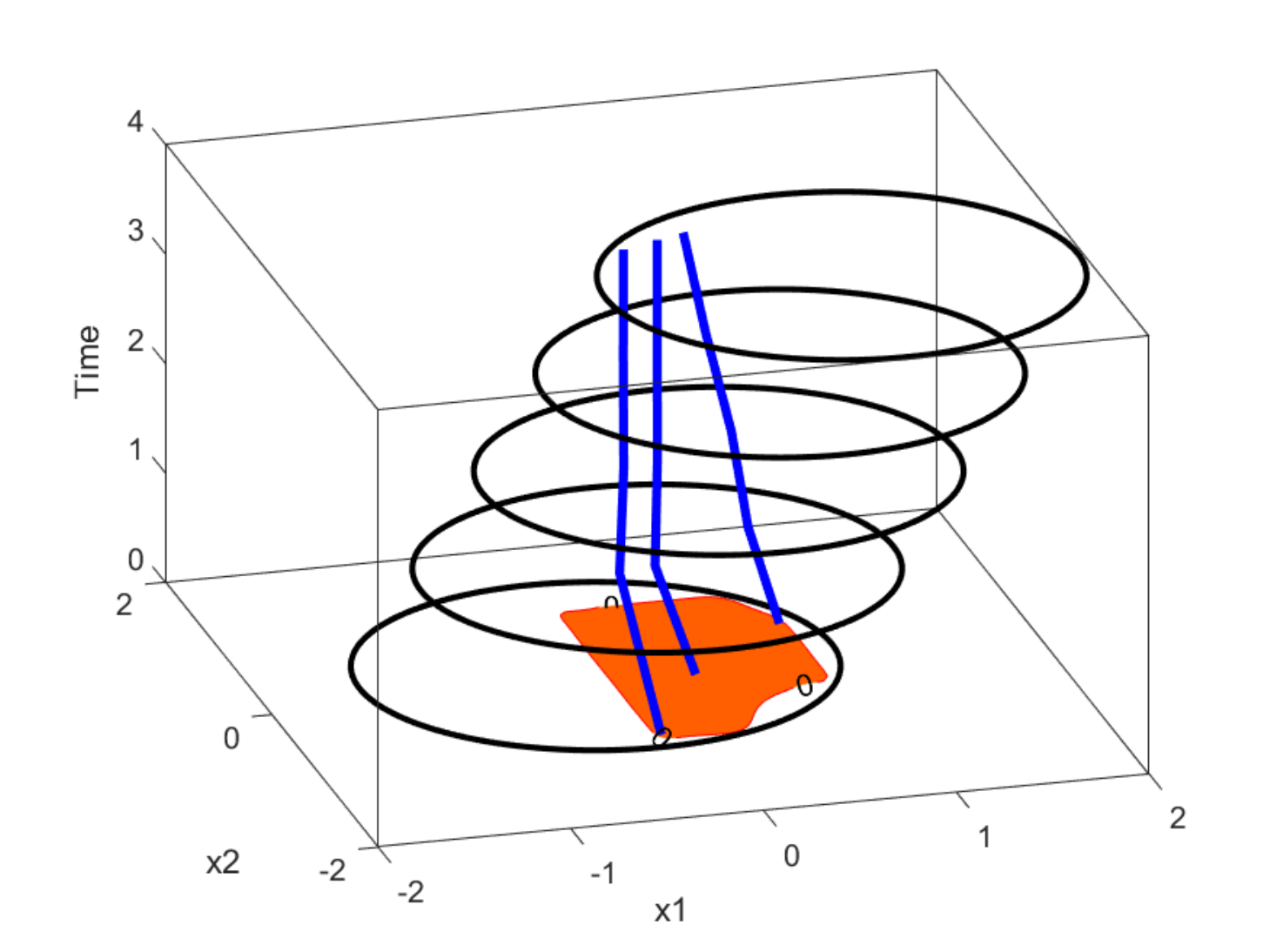}
	\vspace{-20pt}
	\caption{Figure showing an approximation of $L(V,0):=\{x \in \R^n: V(x,0) \le 0\}$, shown in the shaded orange region, where $V$ is the value function of the MSOP associated with Eqn.~\eqref{nex: inv 1}. The z-axis represents time and the black circular lines represent the boundary of $\mathcal A_t$ for $t=1,2,3,4$. Three sample trajectories, shown in blue, start in $L(V,0)$ and remain in the sets $\mathcal A_t$ for the time-steps $t=1,2,3,4$; giving numerical evidence that $L(V,0)$ is indeed an approximation of the FTHMIS.   } \label{fig: inv set 1}
	\vspace{-5pt}
\end{figure}

Value functions can characterize FTHMIS's, as shown by Theorem~\ref{thm: FTHMIS characterization}. We now approximate a FTHMIS by computing a value function using a discretization scheme for solving the GBE~\eqref{eqn: generalize Bellman} using $\{\phi_t\}_{t=0}^T$ as given in Example~\ref{ex: sup forward sep}. Let us consider a discrete time switching system, whose Robust Maximal Invariant Set (RMIS) was previously computed in \cite{xue2018robust}:
\begin{align} \label{nex: inv 1}
f(x,u,t) & = \begin{cases}
\begin{bmatrix}
x_1 \\
(0.5 + u)x_1 - 0.1 x_2
\end{bmatrix} \text{ if } 1-(x_1-1)^2 -x_2^2 \le 0 \\
\begin{bmatrix}
x_2 \\
0.2 x_1 - (0.1 + u)x_2 +x_2^2
\end{bmatrix} \text{ otherwise}.
\end{cases}
\end{align}
We now compute the FTHMIS, denoted by $\mathcal R$, associated with
\begin{align*}
\mathcal A_t & = \{x \in \R^2: g_t(x) \le 0\} \text{ for all } t\in\{0,..,T\},\\
g_t(x) & = \left(x_1 - \frac{(t-1)}{4}\right)^2 + \left(x_2 - \frac{(t+1)}{4}\right)^2 - 1.5,\\
X_t & = [-1,1]^2  \text{ for all } t\in\{0,..,T\},\\
U & = \{u \in \R : u^2 -0.01 \le 0\}, \quad T  = 4.
\end{align*}
Figure~\ref{fig: inv set 1} shows the FTHMIS, $\mathcal R$, found by using a discretization scheme to solve the GBE~\eqref{eqn: generalize Bellman} for $5 \times 5$ state grid points in $[-1,1]^2$. To represent $\mathcal R$ in $\R^2$, once the value function, $V$, is found at each grid point a polynomial function is fitted and its zero-sublevel set, shown as the orange shaded region, approximately gives $\mathcal R$.

\vspace{-0.35cm}
\section{Conclusion} \label{sec: con} 
\vspace{-0.1cm}
For MSOP's with monotonically backward separable cost functions we have derived necessary and sufficient conditions for solutions to be optimal. We have shown that by solving the Generalized Bellman's Equation (GBE) one can derive an optimal input sequence. Furthermore, we have demonstrated the GBE can be numerically solved using a discretization scheme and Approximate Dynamic Programing (ADP) techniques such as Rollout. We have shown our numerical methods can solve current practical problems of interest; such as path planning and the computation of maximal invariant sets. 

\vspace{-0.25cm}

\bibliographystyle{plain}        
\bibliography{bibo}

\parpic{\includegraphics[width=1in,clip,keepaspectratio]{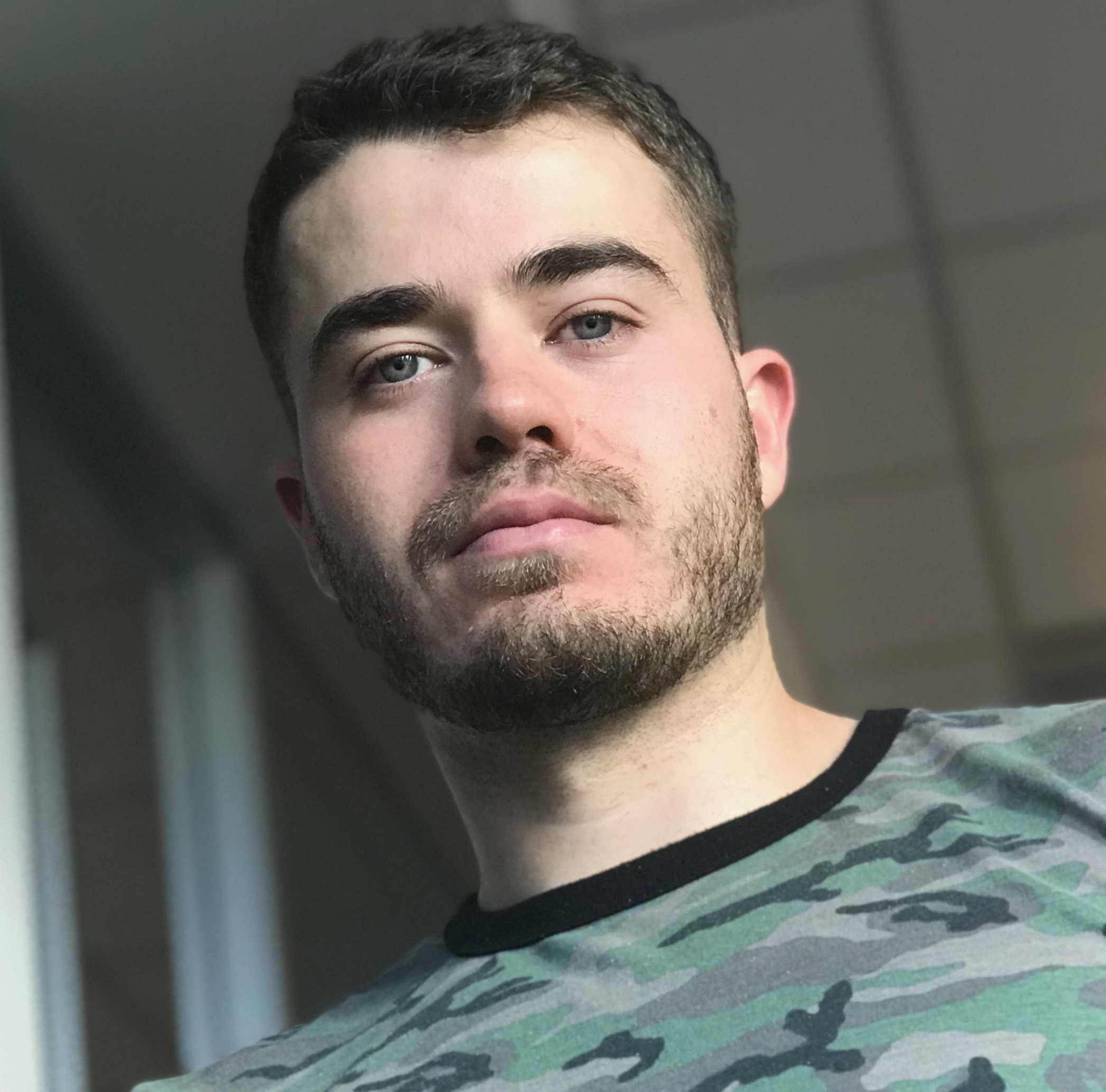}} \small
\noindent {\bf Morgan Jones} received the B.S. and Mmath in
	mathematics from The University of Oxford, England in 2016.
	He is a research associate with Cybernetic
	Systems and Controls Laboratory (CSCL) at ASU. His research primarily focuses on the analysis of nonlinear ODE's and Dynamic Programming.

	\parpic{\includegraphics[width=1in,clip,keepaspectratio]{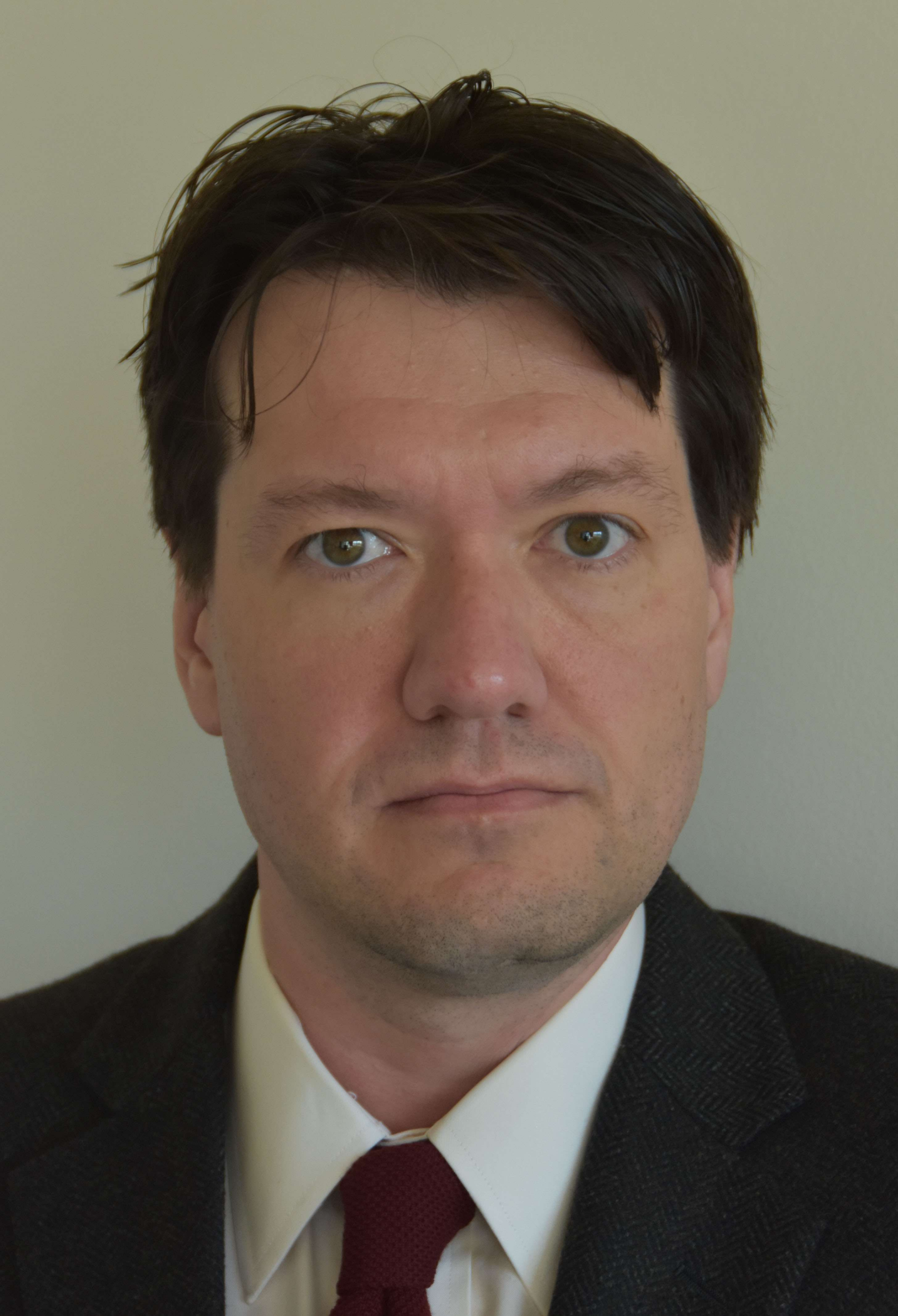}} \small
	\noindent {\bf Matthew M. Peet} received the B.S. degree in
		physics and in aerospace engineering from the University
		of Texas, Austin, TX, USA, in 1999 and
		the M.S. and Ph.D. degrees in aeronautics and astronautics
		from Stanford University, Stanford, CA,
		in 2001 and 2006, respectively. He was a Postdoctoral
		Fellow at INRIA, Paris, France from 2006 to
		2008. He was an Assistant Professor of Aerospace
		Engineering at the Illinois Institute of Technology,
		Chicago, IL, USA, from 2008 to 2012. Currently, he
		is an Associate Professor of Aerospace Engineering
		at Arizona State University, Tempe, AZ, USA. Dr. Peet received a National Science Foundation CAREER award in 2011.

\end{document}